\documentclass[11pt]{article}
\usepackage[T1]{fontenc}
\usepackage{lmodern}
\usepackage[margin=0.9in]{geometry}
\usepackage{amsmath,amssymb,amsthm,mathtools,mathrsfs,bm}
\usepackage{microtype,booktabs,longtable,array,enumitem}
\usepackage{graphicx,xcolor,tikz,tikz-cd}
\usetikzlibrary{arrows.meta,positioning,calc}
\usepackage{needspace,fancyhdr,etoolbox}
\usepackage[colorlinks=true,linkcolor=blue!40!black,citecolor=blue!40!black,urlcolor=blue!40!black,bookmarksnumbered=true]{hyperref}
\hypersetup{pdftitle={Zeta Space: Equivariant Seams, Spectral States, and Fourier--Mellin Duality},pdfsubject={A reference formulation of the zeta-space program}}
\setlist{itemsep=0.25em,topsep=0.4em,leftmargin=1.7em}
\newtheorem{theorem}{Theorem}[section]
\newtheorem{proposition}[theorem]{Proposition}
\newtheorem{lemma}[theorem]{Lemma}
\newtheorem{corollary}[theorem]{Corollary}
\theoremstyle{definition}
\newtheorem{definition}[theorem]{Definition}
\newtheorem{example}[theorem]{Example}

\theoremstyle{remark}
\newtheorem{remark}[theorem]{Remark}
\newtheorem*{mainstatement}{The structural picture}
\AtBeginEnvironment{theorem}{\Needspace{8\baselineskip}}
\AtBeginEnvironment{proposition}{\Needspace{7\baselineskip}}
\AtBeginEnvironment{lemma}{\Needspace{6\baselineskip}}
\AtBeginEnvironment{corollary}{\Needspace{6\baselineskip}}
\numberwithin{equation}{section}
\newcommand{\R}{\mathbb R}\newcommand{\C}{\mathbb C}\newcommand{\Z}{\mathbb Z}
\newcommand{\N}{\mathbb N}\newcommand{\Hh}{\mathbb H}\newcommand{\OO}{\mathbb O}
\newcommand{\LL}{\mathbb L}\newcommand{\FS}{\mathrm{FS}}
\newcommand{\cI}{\mathcal I}\newcommand{\cZ}{\mathcal Z}\newcommand{\cP}{\mathcal P}
\newcommand{\cH}{\mathcal H}\newcommand{\cR}{\mathcal R}
\newcommand{\cD}{\mathcal D}\newcommand{\cB}{\mathcal B}\newcommand{\cM}{\mathcal M}
\newcommand{\Sch}{\mathscr S}\newcommand{\Mer}{\mathscr M}
\newcommand{\dd}{\,d}
\DeclareMathOperator{\Tr}{Tr}\DeclareMathOperator{\Res}{Res}
\DeclareMathOperator{\Hom}{Hom}
\DeclareMathOperator{\PSL}{PSL}
\DeclareMathOperator{\divisor}{div}
\newcommand{\purpose}[1]{\noindent\emph{Purpose.} #1\par\medskip}
\newcommand{\G}{{\Pi_3}}
\newcommand{\Lam}{\Lambda_{\zeta}}

\begin{document}
\hypersetup{pageanchor=false}
\begin{titlepage}
\centering
\vspace*{1.1cm}
{\Huge\bfseries Foundations of Zeta Space: Geometry, Spectral States, and Fourier–Mellin Duality\par}
\vspace{0.55cm}
{\Large John McManus Zimmerman\\[0.25em]johnzimmerman25@gmail.com \par}

\vspace{1cm}
\begin{minipage}{0.91\textwidth}
\begin{abstract}
The zeta-space program seeks a geometric organization of theta states and completed zeta functions in which local spectral data, gluing, and Fourier duality belong to a single compatible construction. Its distinguished three-seed consists of four sheets arranged around an octahedral seam, with eight marked ends and full octahedral symmetry. A cubical boundary presents this incidence structure through six face crosses and eight corner regions.

This manuscript gives precise definitions and proofs for the principal mechanisms of the framework. An equilateral Kirchhoff realization of the seam has an invariant Laplacian unitarily equivalent to a Neumann interval. Its normalized heat trace recovers the Jacobi theta function, and its Mellin transform gives the completed Riemann zeta function. On the sheets, closed energy forms produce conservative seam matching, while hyperbolic cusp collars yield Fourier--Bessel propagation. Boundary transfer, harmonic extension, and ambient hyperbolic Poisson transforms are formulated with their domains and coupling data specified.

For the higher tower, we define enriched bonding morphisms, compatible theta families, and Mellin-compatible transitions. A duality-and-assembly theorem carries Fourier reciprocity to a functional equation over the doubled inverse-limit base. A split rank-two calculation realizes a rank-lowering transition by residues and explains how duality exchanges its two polar prescriptions. The framework thereby separates geometric incidence, analytic propagation, arithmetic specialization, and global assembly while making their relationships explicit.
\end{abstract}
\end{minipage}
\vfill
\begin{minipage}{0.91\textwidth}
\small
\textbf{Reading guide.} The introduction states the program's objectives and the main results. Sections~\ref{sec:geometry}--\ref{sec:rankone} give a complete route from the boundary cube to completed zeta. Sections~\ref{sec:sheetdata}--\ref{sec:maassmellin} develop the sheet and cusp analysis. Sections~\ref{sec:higher}--\ref{sec:assembly} formulate the tower and prove its compatibility results. The final section and appendices record enrichment, conventions, and exact reference formulas.
\end{minipage}
\end{titlepage}
\hypersetup{pageanchor=true}
\pagenumbering{roman}
\tableofcontents
\clearpage
\pagenumbering{arabic}

\section{The aim of the zeta-space program}\label{sec:aim}

A vibrating circle has integer Fourier modes. The heat equation weights those modes by Gaussians. Poisson summation relates short and long scales, and a Mellin transform turns that relation into the functional equation of a completed zeta function. This familiar chain is the analytic prototype of the zeta-space program:
\begin{equation}\label{eq:prototype}
\text{periodic geometry}\ \longrightarrow\ \text{theta state}
\ \longrightarrow\ \text{Mellin state}.
\end{equation}
The program asks how this chain can be organized when the geometry consists of several sheets meeting along a seam, when a finite symmetry permutes the local channels, and when compatible geometries occur at successive ranks.

The central object at the first nontrivial level is the three-seed
\begin{equation}
\mathfrak S_3=(\cI_3,\Gamma_3,\Pi_3).
\end{equation}
Here $\cI_3$ is the sheeted space, $\Gamma_3$ is its seam, and $\Pi_3$ is the symmetry acting on the complete datum. The seam is octahedral. Four distinguished sheets are indexed by the four pairs of opposite cube vertices, and their eight marked ends are indexed by the individual vertices. The cube is therefore a boundary presentation of the incidence geometry.

The objective is to make a spectral state remember this geometry. A local wave should carry information about its sheet and its seam traces. A global state should be assembled from compatible local states. Fourier duality should respect that assembly, and the Mellin transform should carry it to a corresponding symmetry of completed states. At higher levels, these structures should persist under explicitly chosen bonding morphisms.

The spectral interpretation of zeta functions has important precedents
in Connes's noncommutative trace formula and Deninger's proposed
connections between arithmetic and dynamics on foliated spaces
\cite{Connes1999,Deninger1998}.
These provide context for the program's spectral objectives;
the constructions developed here begin with the specified seed
geometry and its operators.

\subsection{Three objectives and their mathematical outputs}

The first objective is \emph{geometric}: specify the seed, its symmetry, and its incidence maps. Its output is a stratified space with enough metric and bundle data to define operators. The second is \emph{analytic}: identify which local waves glue and how their Fourier or theta descriptions transform under scale. Its output is a space of compatible spectral states and a Mellin construction. The third is \emph{global}: assemble the constructions over a tower and make duality act on the assembled object. Its output is the organization
\begin{equation}\label{eq:threelevels-intro}
\cZ\longrightarrow\Theta\longrightarrow\LL,
\end{equation}
where $\LL$ is the geometric base, $\Theta$ records theta data, and $\cZ$ records those data together with their completed Mellin states.

The arithmetic aim is to study completed zeta functions, their functional symmetries, and their divisors through this organization. The rank-one reference state is the Riemann zeta function. Its eventual spectral-localization target is the Riemann hypothesis: an identity that locates the nontrivial divisor on the critical line. Proposition~\ref{prop:detcriterion} formulates one precise criterion for such a localization. Higher-rank states require their own lattices, kernels, measures, and transition laws. The framework treats these choices as mathematical structure, so that a statement about a zeta function can be traced back to the data that produce it.

\Needspace{12\baselineskip}
\subsection{The main results and their scope}

\begin{mainstatement}
The reference formulation has three complementary components.
\begin{enumerate}
\item For the equilateral octahedral seam with Kirchhoff conditions, its $\Pi_3$-invariant Laplacian is a Neumann interval operator. A normalized heat trace is the Jacobi theta function, and its Mellin transform is $\pi^{-s/2}\Gamma(s/2)\zeta(s)$; see Theorems~\ref{thm:interval} and~\ref{thm:seamzeta}.
\item For the specified sheetwise energy form, a unique nonnegative self-adjoint operator governs propagation. Its weak gluing law is continuity of trace together with conservation of conormal flux; see Theorem~\ref{thm:sheetoperator}. Cusp collars and global automorphic states have the Fourier--Bessel and Mellin descriptions of Theorems~\ref{thm:cusp} and~\ref{thm:maassFE}.
\item An enriched tower whose state transitions commute with Mellin transformation and whose duality respects those transitions has an assembled completed state satisfying the induced reflection law; see Theorem~\ref{thm:assembly}. The rank-two residue calculation in Theorem~\ref{thm:ranktwo} gives an explicit model of such a transition.
\end{enumerate}
\end{mainstatement}

Each component has its own operator and domain. The seam operator is an exactly solvable model attached to $\Gamma_3$. The sheet operator depends on the sheet metrics and their transmission conditions. The tower theorem applies to enriched systems satisfying its stated compatibility identities. Keeping these inputs visible makes the manuscript usable as a reference: every result can be applied by checking a finite list of structural hypotheses.

\subsection{Familiar mathematics and the organizing principle}

The analytic tools are classical: metric-graph Laplacians, closed quadratic forms, Fourier series, Poisson summation, hyperbolic cusp expansions, and Mellin transforms. Quantum-graph theory provides the natural operator language for a metric seam~\cite{Kuchment}; the theta--Mellin mechanism is the classical prototype behind Riemann's and Tate's functional equations~\cite{Poonen}. The role of zeta space is to organize these mechanisms over a common geometric datum and to state the conditions under which they pass through gluing and change of rank.

Three recurring distinctions keep that organization precise. An incidence identification specifies which points meet; a transmission law specifies how fields meet. Equivariance specifies how a group acts on data; invariance selects states fixed by that action. Finally, the parameter of a Laplace eigenvalue and the variable of a Mellin transform serve different purposes. We use $\eta=1/2+ir$ for the former and $s$ for the latter. The notation $\LL$ is reserved for the geometric base throughout.

\section{The cube, the seam, and the distinguished seed}\label{sec:geometry}
\purpose{Recover the finite incidence structure that every later analytic construction must respect.}

\subsection{Opposite corners and the full foliated system}

Let $X=[0,1]^3$, let $V$ be its set of eight vertices, and write $v^\ast$ for the vertex opposite $v$. Thus $v$ and $v^\ast$ are the two endpoints of a body diagonal. Define
\begin{equation}
\cP=\{\{v,v^\ast\}:v\in V\}.
\end{equation}
This set contains four pairs. The symbols $v$ and $v^\ast$ are geometric labels; no coordinate representation is needed at this stage.

\begin{definition}[Foliated incidence datum]
For each $p\in\cP$, let $\mathcal A_p$ index a family of open leaves $L_{p,a}\cong(0,1)\times S^1$, with their two ends labeled by the two vertices in $p$. An incidence relation on their end compactifications defines
\begin{equation}\label{eq:FS}
\FS(X)=\left(\bigsqcup_{p\in\cP}\ \bigsqcup_{a\in\mathcal A_p}
\widehat L_{p,a}\right)\Big/\!\sim_{\mathrm{inc}}.
\end{equation}
A realization includes the topology, the incidence maps, and a $\Pi_3$ action preserving this structure. Its regular regions may be foliated; its shared strata are recorded explicitly by $\sim_{\mathrm{inc}}$.
\end{definition}

The three-seed selects a distinguished equivariant four-sheeted subsystem inside this larger object. Equation~\eqref{eq:FS} retains the full family: a calculation near one cusp describes one local analytic channel within it.

\subsection{The six crosses}

On each cube face $F$, choose coordinates $(u,w)\in[0,1]^2$ and set
\begin{equation}
C_F=\{u=1/2\}\cup\{w=1/2\},
\qquad C_{\partial X}=\bigcup_F C_F.
\end{equation}
The four arms of $C_F$ join its center to its edge midpoints. Arms on adjacent faces meet at the shared midpoint.

The topological constructions used below---including graph homology,
attaching maps, and quotient spaces---are understood in the standard
sense of algebraic topology \cite{Hatcher2002}.

\begin{proposition}[Boundary realization of the seam]\label{prop:cross}
The graph $C_{\partial X}$ is a subdivision of the octahedral graph $K_{2,2,2}$. It has six degree-four vertices, twelve degree-two vertices, and twenty-four subdivided edges. Suppressing the degree-two vertices gives a graph with six vertices and twelve edges, and
\begin{equation}
\overline C_{\partial X}\cong\Gamma_3:=K_{2,2,2},
\qquad b_1(\Gamma_3)=7.
\end{equation}
These identifications are equivariant for the cube's full symmetry group $\Pi_3\cong B_3\cong O_h$.
\end{proposition}
\begin{proof}
Each face center has four incident arms. Each of the twelve cube edges contributes one common midpoint, of degree two. Removing the subdivision leaves one vertex for each face and one edge for each pair of adjacent faces. A face is adjacent to every face except itself and its opposite. The three pairs of opposite faces are therefore the three parts of $K_{2,2,2}$. The graph is connected, so its first Betti number is $12-6+1=7$. Every construction uses face adjacency and midpoints, both preserved by cube symmetries. Finally, those symmetries are the signed permutations of three centered coordinate axes, a group of order $2^3\cdot3!=48$.
\end{proof}

\begin{remark}[Linearization preserves the graph]
Suppressing a subdivision is a change of cell structure. It induces a homeomorphism of the geometric realizations, and an edge parametrization can identify their lengths. Thus calling $C_{\partial X}$ a linearized seam means that it displays the seam in planar face charts; it does not collapse its cycles.
\end{remark}

\begin{proposition}[Eight triangular channels]\label{prop:triangles}
The complement of $C_{\partial X}$ in the cube boundary has eight connected regions, naturally indexed by $V$. Their boundaries correspond to the eight triangular face cycles of the octahedron. Every seam edge occurs in exactly two such cycles. Opposite cube vertices determine four pairs of opposite triangular regions.
\end{proposition}
\begin{proof}
Each face cross cuts its face into four quadrants. The three quadrants incident to a given cube vertex meet across the original cube edges, forming one region around that vertex. Midpoints separate this region from the analogous regions at neighboring vertices. Hence the $24$ quadrants form $8$ regions of $3$ quadrants each. If the faces $F_1,F_2,F_3$ meet at $v$, the boundary becomes the cycle through $c_{F_1},c_{F_2},c_{F_3}$ in the reduced graph. A cube edge has two endpoints, giving the two triangles containing its dual seam edge. The antipodal correspondence follows from the cube's central inversion.
\end{proof}

\begin{corollary}[Metric triangular cycles]\label{cor:trianglelength}
If every seam edge has length $\ell>0$, each of the eight triangular cycles has length $3\ell$. The girth of the equilateral metric graph, measured as the shortest embedded cycle length, is $3\ell$.
\end{corollary}
\begin{proof}
Every triangle uses three distinct edges. The graph is simple and has no loops or double edges, so no cycle uses fewer than three edges.
\end{proof}

\begin{center}
\renewcommand{\arraystretch}{1.14}
\begin{tabular}{@{}lll@{}}
\toprule
Boundary datum & Seam datum & Spectral role\\
\midrule
Six faces & Six vertices & Junction labels\\
Twelve cube edges & Twelve seam edges & Transmission channels\\
Eight corners & Eight triangular regions & Distinguished ends\\
Four opposite pairs & Four opposite region pairs & Distinguished sheets\\
Twenty-four quadrants & Three incidences per region & Boundary inputs\\
\bottomrule
\end{tabular}
\end{center}

\subsection{The four-sheeted object}

\begin{definition}[Three-seed with specified gluing]\label{def:seed}
A three-seed is a datum $\mathfrak S_3=(\cI_3,\Gamma_3,\Pi_3)$ with
\begin{equation}
\cI_3=\left(\bigsqcup_{p\in\cP}X_p\right)\Big/\!\sim_{\Gamma_3},
\qquad \Gamma_3\cong K_{2,2,2},
\qquad \Pi_3\cong O_h.
\end{equation}
The datum includes the marked seam strata on the sheets, their attaching maps, and an action preserving the resulting incidence structure. Each sheet has two distinguished ends, labeled by the pair $p$. Conical compactifications and cusp metrics are analytic enrichments specified below.
\end{definition}

The quotient notation records attaching maps. A quotient by a group action is a different presentation and may be used when that action has been supplied and its orbit relation agrees with the intended identifications. Equivariance alone specifies how symmetries act on the seed; it does not replace the attaching maps.

\begin{figure}[t]
\centering
\includegraphics[width=0.47\textwidth]{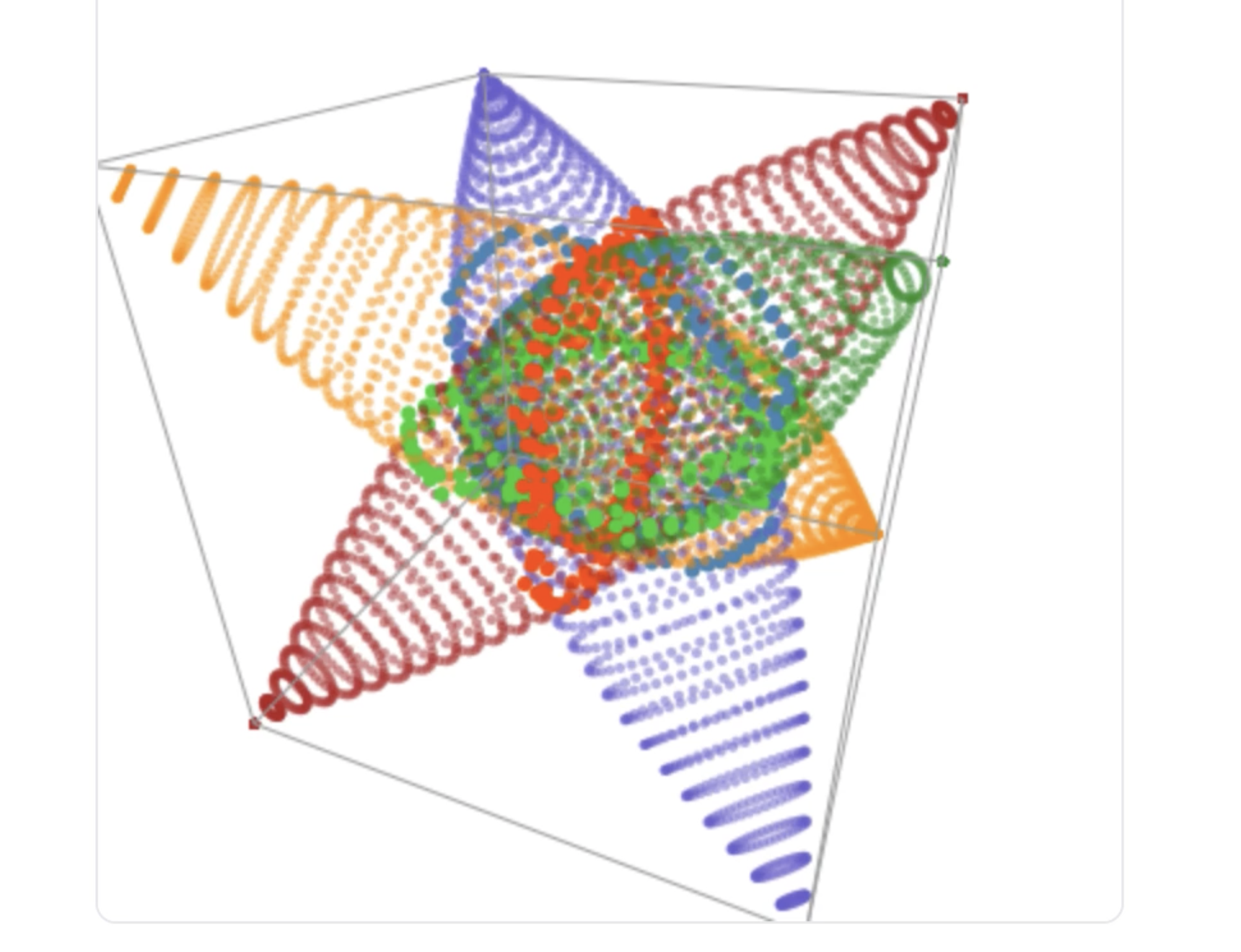}\hfill
\includegraphics[width=0.44\textwidth]{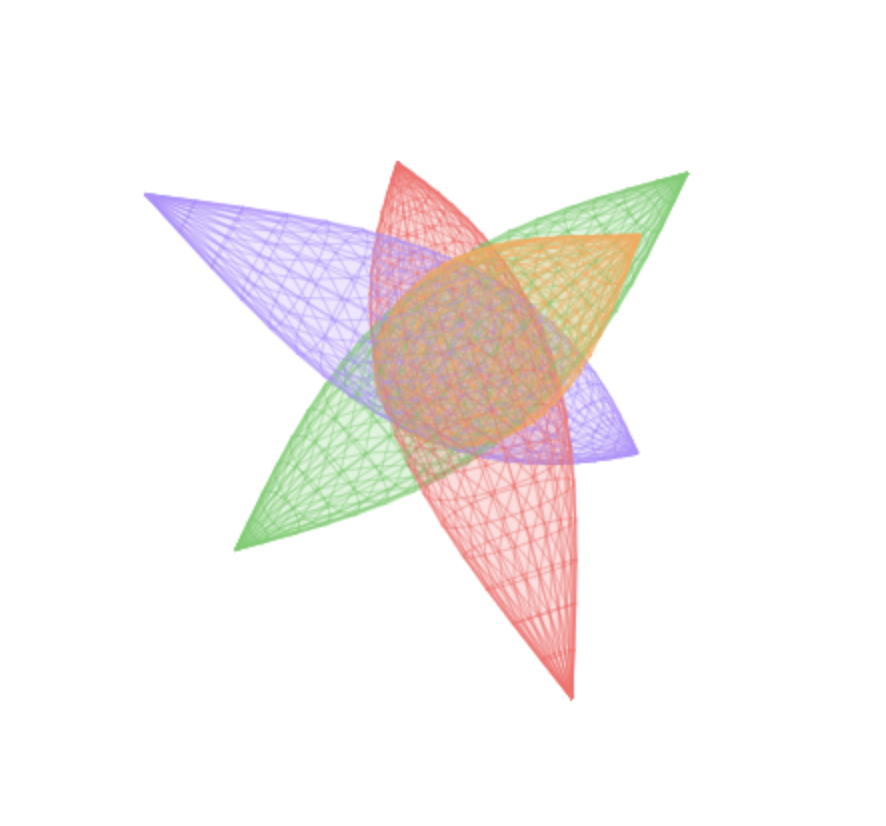}
\caption{Provided visualizations of circular leaf families and a four-sheeted arrangement. They guide the geometric intuition for $\FS(X)$ and its distinguished subsystem. The seam incidence, intrinsic metrics, and operator domains used in the proofs are specified in the text.}
\label{fig:seedpictures}
\end{figure}

\section{An exactly solvable invariant seam}\label{sec:seamoperator}
\purpose{Attach a concrete self-adjoint operator to the seam and identify the spectral channel that will produce the rank-one theta state.}

\subsection{The equilateral Kirchhoff operator}

We use the standard framework of quantum graphs: each edge carries
a differential operator, and vertex conditions specify its global
self-adjoint realization
\cite{BerkolaikoKuchment2013,KostrykinSchrader1999}.

Give each edge of $\Gamma_3$ length $\ell$. Let
\begin{equation}
H_\Gamma=\bigoplus_{e\in E(\Gamma_3)}L^2(0,\ell),
\qquad
V_\Gamma=\{f\in\bigoplus_eH^1(0,\ell):f\text{ is continuous at vertices}\}.
\end{equation}
The energy is $q_\Gamma(f)=\sum_e\int_0^\ell|f'_e|^2\dd x$.

\begin{lemma}[Closed seam energy]\label{lem:graphform}
The form $q_\Gamma$ on $V_\Gamma$ is densely defined, nonnegative, and closed. Its associated operator $A_\Gamma$ acts as $-d^2/dx^2$ on each edge, with domain consisting of edgewise $H^2$ functions satisfying continuity and
\begin{equation}\label{eq:Kirchhoff}
\sum_{e\ni v}\partial_e f(v)=0
\end{equation}
at every vertex. Here $\partial_e$ points from the vertex into the edge. The operator has compact resolvent.
\end{lemma}
\begin{proof}
The form norm is the direct-sum $H^1$ norm. Vertex evaluation is continuous on each $H^1$ interval, so the continuity constraints define a closed subspace. Functions compactly supported in edge interiors are dense in the direct-sum $L^2$ space, proving density. The representation theorem for closed nonnegative forms gives a unique nonnegative self-adjoint operator. Testing against interior functions gives $-f_e''$ and edgewise $H^2$ regularity. Integration by parts leaves at each vertex a boundary coefficient proportional to the sum in~\eqref{eq:Kirchhoff}; the common test value at that vertex is arbitrary, so that coefficient vanishes. Conversely, those conditions remove every boundary term. Finally, a finite direct sum of $H^1$ intervals embeds compactly into the direct-sum $L^2$ space. The form-domain embedding, and hence the resolvent, is compact.
\end{proof}

The group $\G$ acts unitarily by pullback, preserves $V_\Gamma$, and preserves the energy. Hence its Reynolds projection
\begin{equation}
P_\G=\frac1{48}\sum_{g\in\G}U_g
\end{equation}
commutes with $A_\Gamma$ and its heat semigroup. Denote the restriction to $H_\Gamma^\G$ by $A_\Gamma^\G$.

Averaging over the finite symmetry group gives the orthogonal
projection onto invariant vectors, as in the standard representation
theory of finite groups \cite{Serre1977}.

\begin{theorem}[Invariant seam reduction]\label{thm:interval}
The operator $A_\Gamma^\G$ is unitarily equivalent to the Neumann Laplacian on $[0,\ell/2]$. Its eigenvalues are simple and equal to
\begin{equation}\label{eq:invariantspectrum}
\lambda_m=\left(\frac{2\pi m}{\ell}\right)^2,
\qquad m=0,1,2,\ldots.
\end{equation}
On every edge, the corresponding invariant eigenfunction is proportional to $\cos(2\pi mx/\ell)$.
\end{theorem}
\begin{proof}
Signed permutations act transitively on oriented edges. There is also a symmetry that reverses a chosen edge: for the edge from $e_1$ to $e_2$, exchange the first two coordinates. Therefore an invariant function is represented on every edge by the same profile $F$, independent of orientation, and $F(x)=F(\ell-x)$. Conversely, any such profile gives an invariant edgewise function.

At the form level, continuity holds automatically because $F(0)=F(\ell)$. Restriction to $[0,\ell/2]$ determines $F$ by reflection, and
\begin{equation}
\|f\|^2_{H_\Gamma}=24\int_0^{\ell/2}|F(x)|^2\dd x,
\qquad q_\Gamma(f)=24\int_0^{\ell/2}|F'(x)|^2\dd x.
\end{equation}
Thus $f\mapsto\sqrt{24}\,F|_{[0,\ell/2]}$ is a unitary equivalence of forms with the unrestricted $H^1$ energy on that interval. Its associated operator is the Neumann Laplacian. Its eigenfunctions are $\cos(\pi mx/(\ell/2))$, giving~\eqref{eq:invariantspectrum}. Equivalently, the Kirchhoff condition gives $4F'(0)=0$, while reflection gives $F'(\ell/2)=0$.
\end{proof}

\begin{remark}[The origin of rank one]
The graph has seven independent cycles, but its invariant scalar operator has one interval coordinate. The theta rank in this model comes from the integer mode label of the symmetry-reduced operator. It is not the first Betti number of the full seam.
\end{remark}

\begin{remark}[A general symmetry reduction]
The same proof applies to any finite equilateral metric graph with a group acting transitively on its directed edges. The norm factor becomes twice the number of edges, and the invariant operator is again the Neumann interval of length $\ell/2$. The octahedral seed places this common spectral channel inside its particular four-sheeted, eight-ended incidence structure.
\end{remark}

\begin{corollary}[Heat trace and normalization]\label{cor:seamheat}
Set $\mathsf H_\Gamma=\ell^2A_\Gamma^\G/(4\pi)$ and define
\begin{equation}\label{eq:seamtheta}
Z_\Gamma(t)=\Tr(e^{-t\mathsf H_\Gamma}),
\qquad \theta_\Gamma(t)=2Z_\Gamma(t)-1.
\end{equation}
Then, for $t>0$,
\begin{equation}
Z_\Gamma(t)=\sum_{m=0}^\infty e^{-\pi m^2t},
\qquad
\theta_\Gamma(t)=\sum_{m\in\Z}e^{-\pi m^2t}.
\end{equation}
\end{corollary}
\begin{proof}
The normalized eigenvalues are $\pi m^2$. The heat trace is their convergent exponential sum. Doubling it and subtracting one restores the two signs of each nonzero integer while counting the constant mode once.
\end{proof}

This gives a fully specified version of the rank-one extraction. The chosen object is the invariant heat trace of the Kirchhoff seam, with a stated normalization. It is distinct from either the full graph heat trace or an unspecified sheet heat trace.

Heat kernels on metric graphs also admit descriptions in terms of
walks and associated trace formulas
\cite{KostrykinPotthoffSchrader2007}.
Here the symmetry reduction permits a direct calculation of the
invariant heat trace.

\begin{proposition}[The triangular period and its harmonics]\label{prop:triangleharmonics}
Restriction of an invariant seam eigenfunction to any triangular face boundary is periodic with period $\ell$ along a circle of length $3\ell$. In the full Fourier basis of that circle, only indices divisible by $3$ occur.
\end{proposition}
\begin{proof}
Each edge carries $F(x)=\cos(2\pi mx/\ell)$, and reversal does not alter it. Concatenating the three edge profiles therefore gives the same $\ell$-periodic function three times. Its circular frequency $2\pi m/\ell$ equals $2\pi(3m)/(3\ell)$.
\end{proof}

\begin{remark}[Length factors]
A length scale contributes a Mellin factor only after the operator and spectral variable are fixed. For the positive spectrum above,
\begin{equation}
\Tr{}'\bigl((A_\Gamma^\G)^{-z}\bigr)
=\left(\frac{\ell}{2\pi}\right)^{2z}\zeta(2z),\qquad \Re z>\tfrac12.
\end{equation}
The triangular length $3\ell$ and the seam period $\ell$ are both geometrically meaningful; Proposition~\ref{prop:triangleharmonics} explains their relation before any factor such as $3^{-s}$ is introduced.
\end{remark}

\section{From the invariant seam to completed zeta}\label{sec:rankone}
\purpose{Complete the first objective: derive the theta reciprocity and Mellin functional equation for the explicit seam state.}

The classical theory of theta functions supplies the transformation
laws underlying the calculation \cite{Mumford2007}.
Its interpretation through Fourier analysis and zeta integrals belongs
to the broader framework of Tate's thesis \cite{Kudla2004}.

\begin{lemma}[Gaussian Fourier transform and Poisson summation]\label{lem:poisson}
For the Fourier convention $\widehat f(\xi)=\int_\R f(x)e^{-2\pi ix\xi}\dd x$,
\begin{equation}
\widehat{e^{-\pi t x^2}}(\xi)=t^{-1/2}e^{-\pi\xi^2/t},\qquad t>0.
\end{equation}
For a Schwartz function $f$,
$\sum_{m\in\Z}f(m)=\sum_{m\in\Z}\widehat f(m)$.
\end{lemma}
\begin{proof}
At $t=1$, differentiate the Fourier integral in $\xi$ and integrate by parts in $x$. The result satisfies $F'(\xi)=-2\pi\xi F(\xi)$, while $F(0)=\int e^{-\pi x^2}\dd x=1$, as follows from the two-dimensional Gaussian integral in polar coordinates. Thus $F(\xi)=e^{-\pi\xi^2}$; scaling proves the formula for $t>0$. For Poisson summation, periodize $f$ as $\sum_m f(x+m)$. Schwartz decay permits termwise differentiation and integration; its Fourier coefficient at $n$ is $\widehat f(n)$. The Fourier series converges absolutely, and evaluation at $x=0$ gives the identity.
\end{proof}

\begin{theorem}[Reciprocity of the seam theta state]\label{thm:thetareciprocity}
The state in~\eqref{eq:seamtheta} satisfies
\begin{equation}\label{eq:thetareciprocity}
\theta_\Gamma(t)=t^{-1/2}\theta_\Gamma(1/t).
\end{equation}
\end{theorem}
\begin{proof}
Apply Lemma~\ref{lem:poisson} to $e^{-\pi t x^2}$ and use Corollary~\ref{cor:seamheat}.
\end{proof}

\begin{theorem}[Completed zeta from the invariant seam]\label{thm:seamzeta}
For $\Re s>1$, the Mellin state
\begin{equation}\label{eq:seamMellin}
\Lam(s):=\int_0^\infty\bigl(Z_\Gamma(t)-1\bigr)t^{s/2}\frac{\dd t}{t}
\end{equation}
is equal to $\pi^{-s/2}\Gamma(s/2)\zeta(s)$. It has a meromorphic continuation to $\C$ with simple poles at $0,1$, of residues $-1,1$, and satisfies
\begin{equation}\label{eq:zetaFE}
\Lam(s)=\Lam(1-s).
\end{equation}
The function $\xi(s)=\tfrac12s(s-1)\Lam(s)$ is entire and obeys $\xi(s)=\xi(1-s)$.
\end{theorem}
\begin{proof}
For real $s>1$, Tonelli's theorem gives
\begin{align}
\int_0^\infty\sum_{m\ge1}e^{-\pi m^2t}t^{s/2}\frac{\dd t}{t}
&=\sum_{m\ge1}(\pi m^2)^{-s/2}\Gamma(s/2)\\
&=\pi^{-s/2}\Gamma(s/2)\zeta(s).
\end{align}
The same calculation bounds the absolute value for complex $s$ by the expression at $\Re s$, proving absolute convergence and the identity on $\Re s>1$.

Using $Z_\Gamma-1=(\theta_\Gamma-1)/2$, split the integral at $1$. On $(0,1)$ insert~\eqref{eq:thetareciprocity} and substitute $t\mapsto1/t$ in the decaying part. This yields
\begin{equation}\label{eq:splitMellin}
\Lam(s)=\frac1{s-1}-\frac1s+
\frac12\int_1^\infty(\theta_\Gamma(t)-1)
\left(t^{s/2}+t^{(1-s)/2}\right)\frac{\dd t}{t}.
\end{equation}
The integral is entire: on every compact $s$-set its integrand and all $s$-derivatives are bounded by a polynomial in $t$ and $\log t$ times $Ce^{-\pi t}$. The rational part supplies exactly the stated poles and residues. Both parts of~\eqref{eq:splitMellin} are invariant under $s\mapsto1-s$. Multiplication by $s(s-1)/2$ cancels the two poles and preserves the reflection.
\end{proof}

For the classical analytic theory of $\zeta(s)$, its completed
functional equation, and the zeros of $\xi(s)$, we refer to
\cite{Titchmarsh1986}.

\begin{remark}[What the seed contributes]
The theta and zeta identities are classical. Theorem~\ref{thm:interval} identifies a specific equivariant seam operator whose normalized heat trace realizes their input. Thus the rank-one reference state is attached to an explicit analytic enrichment of $\Gamma_3$, rather than inferred from the seam's symmetry alone.
\end{remark}

\begin{proposition}[Eight scalar channels and Reynolds projection]\label{prop:eightchannels}
Let $\G$ act on an eight-component scalar kernel by permuting its vertex labels. The invariant subspace consists of diagonal kernels $(a(t),\ldots,a(t))$, and projection onto it replaces the components by their average. In particular, eight copies of $\theta_\Gamma$ project to that same diagonal state.
\end{proposition}
\begin{proof}
The action on vertices is transitive. Invariance therefore makes all components equal. Averaging over the group gives each original component the same weight $1/8$ in every output component.
\end{proof}

Scalar channel averaging is useful after a kernel has been identified. Gluing spatial wavefunctions additionally requires the trace conditions developed below; a scalar heat trace carries no pointwise edge coordinate.

\section{Analytic enrichments of the four sheets}\label{sec:sheetdata}
\purpose{Specify what must be added to the incidence seed before a differential equation or a spectral statement is meaningful.}

\subsection{Conical compactification and cusp enrichment}

The compact sheet notation $S^2(\alpha,\alpha)$ will denote a sphere with two marked conical points of total angle $\alpha$. When $\alpha=2\pi/m$ with $m\in\N$, this is an orbifold signature in the usual effective sense. For general positive $\alpha$, it is a conical metric signature. The two marked points are the end compactification of a cylinder.

Conical metrics and orbifold structures carry different local data:
the former prescribe metric cone angles, whereas the latter require
local finite-group quotient charts.
Background on these two settings is provided by
\cite{Troyanov1991} and \cite{AdemLeidaRuan2007}, respectively.

\begin{proposition}[An action--angle sheet]\label{prop:actionangle}
Let $a,b,T>0$, let $u\in(0,a/b)$ and $v\in\R/T\Z$, and set
\begin{equation}\label{eq:actionangle}
g_\phi=\frac{du^2}{\phi(u)}+\phi(u)\,dv^2,
\qquad \phi(u)=au-bu^2.
\end{equation}
The metric has area form $du\,dv$, Gaussian curvature $b$, and scalar curvature $2b$. Its metric completion has two conical points, each of angle $aT/2$, and total area $aT/b$.
\end{proposition}
\begin{proof}
The determinant in $(u,v)$ coordinates is one. Let $r$ be meridian arclength, so $dr=du/\sqrt\phi$, and write $f(r)=\sqrt{\phi(u(r))}$. Then $g=dr^2+f(r)^2dv^2$, with $f'=\phi'/2$ and $f''=\phi''f/2$. The curvature is $-f''/f=-\phi''/2=b$. Near $u=0$, set $u=ar^2/4$ at leading order. The metric becomes $dr^2+(a^2/4)r^2dv^2$ to leading order, giving cone angle $aT/2$. The other endpoint has slope $\phi'(a/b)=-a$ and the same angle. Integrating the area form gives $T(a/b)$.
\end{proof}

\begin{remark}[A normalization used in the seed]
For $a=2/\sqrt3$ and angular period $T=2\pi$, the cone angle is $2\pi/\sqrt3$. This is a conical normalization; the word orbifold requires the additional integral-order condition above. The calculation fixes the meaning of the angle rather than leaving it implicit in the symbol $\alpha$.
\end{remark}

The action--angle description used here is a two-dimensional instance
of the symplectic-coordinate approach to toric K\"ahler geometry
\cite{Abreu2003}.

For the Maass realization, the marked ends are equipped with cusp collars
\begin{equation}\label{eq:cuspmetric}
C_Y=(\R/\Z)\times(Y,\infty),\qquad
 g_C=y^{-2}(dx^2+dy^2)=dt^2+e^{-2t}dx^2,\quad t=\log y.
\end{equation}
A compact middle region connects the two ends of each sheet and contains its marked seam strata. The core metric is part of the enrichment and may have variable curvature.

\begin{lemma}[Two ends and the choice of core metric]\label{lem:twocusps}
A complete finite-area surface homeomorphic to a cylinder cannot have constant curvature $-1$ everywhere. A cylindrical sheet can carry two hyperbolic cusp collars joined through a smooth compact region of variable curvature.
\end{lemma}
\begin{proof}
For a complete finite-area hyperbolic surface, Gauss--Bonnet on cusp truncations gives $-\operatorname{area}=2\pi\chi$. The horocyclic boundary contributions tend to zero as the cusps are truncated farther out. A cylinder has Euler characteristic zero, contradicting positive area. For the second statement, write a cylinder as $\R\times S^1$ with metric $dr^2+f(r)^2dx^2$, where $f>0$ is smooth and equals a positive multiple of $e^{-|r|}$ outside a compact interval. The ends are cusp metrics after shifting $r$, and the intervening curvature is $-f''/f$.
\end{proof}

This explains how the conical and Maass pictures are related. They share marked sheet topology and end labels. Their local metrics and spectral operators are different enrichments. In the rest of the sheet analysis, the ends are cusps.

\subsection{Standing analytic datum}

\begin{definition}[Admissible sheet enrichment]\label{def:analyticseed}
An admissible scalar enrichment of the three-seed consists of:
\begin{enumerate}
\item a finite collection of smooth Riemannian surface pieces $M_\alpha$, obtained by cutting the sheets along the seam, with piecewise smooth compact seam boundaries and the prescribed cusp ends;
\item arclength identifications of incident seam edges, with positive constant weights $\kappa_\alpha$;
\item a $\Pi_3$ action carrying the pieces, metrics, weights, and attaching maps to one another;
\item chosen external boundary conditions, if external boundary is present.
\end{enumerate}
The Sobolev spaces are defined using the intrinsic metrics. The seam traces satisfy the usual bounded trace maps into the trace spaces of the compact boundary arcs. At boundary corners, these are understood as the traces of global $H^1$ functions on each piece.
\end{definition}

The trace-space formulation incorporates endpoint compatibility automatically. It avoids assigning independent arbitrary edge traces at a corner when no $H^1$ function realizes them.

\section{Variational gluing and harmonic extension}\label{sec:gluing}
\purpose{Turn the phrase ``local states glue'' into a solvable operator problem and identify the seam-supported matching defect.}

\subsection{The global operator}

For the datum of Definition~\ref{def:analyticseed}, use
\begin{equation}
H_3=\bigoplus_\alpha L^2(M_\alpha,\kappa_\alpha\dd\mathrm{vol}_\alpha).
\end{equation}
Let $V_3$ be the subspace of $\bigoplus_\alpha H^1(M_\alpha)$ whose identified seam traces agree, with homogeneous external Dirichlet conditions where specified. Define
\begin{equation}\label{eq:sheetform}
q_3(u,v)=\sum_\alpha\kappa_\alpha\int_{M_\alpha}
\langle\nabla u_\alpha,\nabla\overline v_\alpha\rangle\dd\mathrm{vol}_\alpha.
\end{equation}

The passage from a densely defined, closed, nonnegative quadratic form
to its associated self-adjoint operator uses the representation
theorem for closed forms \cite{Kato1995}.
The trace and weak boundary-value formulations are based on standard
Sobolev and elliptic theory \cite{McLean2000}.

\begin{theorem}[Conservative sheet operator]\label{thm:sheetoperator}
The form $q_3$ is densely defined, nonnegative, and closed. It determines a unique nonnegative self-adjoint operator $A_3$. On smooth portions of its domain, $A_3u=-\Delta u$ on each piece, and along every open seam edge $e$ its matching conditions are
\begin{equation}\label{eq:sheetmatching}
u_\alpha|_e=h_e\quad(\alpha\in I(e)),\qquad
\sum_{\alpha\in I(e)}\kappa_\alpha\partial_{\nu_\alpha}u_\alpha|_e=0.
\end{equation}
Here $\nu_\alpha$ is the outward conormal of the cut piece. The operator and its functional calculus commute with $\Pi_3$.
\end{theorem}
\begin{proof}
The norm $\|u\|_{H_3}^2+q_3(u,u)$ is equivalent to the finite direct-sum $H^1$ norm, since all weights are fixed and positive. The trace identifications and external homogeneous conditions are kernels of bounded trace maps, so $V_3$ is closed in that norm. Interior compactly supported functions are dense in $H_3$, giving density of the form. The representation theorem gives a unique nonnegative self-adjoint operator.

Interior test functions identify its differential expression. For a sufficiently regular domain element, integration by parts against $v\in V_3$ produces the common edge test trace multiplied by the weighted sum of conormal derivatives. This sum must vanish in the dual trace space. Conversely, the stated differential expression and vanishing boundary terms satisfy the defining weak identity of the form operator. The weak domain supplied by the form also fixes behavior at seam vertices.

Finally, the group action preserves both the Hilbert inner product and the form. If $u=(A_3+1)^{-1}f$, its weak defining equation transforms to the same equation with right side $U_gf$. Uniqueness gives $U_g(A_3+1)^{-1}=(A_3+1)^{-1}U_g$, and the spectral functional calculus yields the remaining commutations.
\end{proof}

\begin{corollary}[Invariant evolution and energy conservation]
The invariant subspace $H_3^{\Pi_3}$ is reducing for $A_3$. The heat evolution $e^{-tA_3}$ and Schr\"odinger evolution $e^{-itA_3}$ preserve it. For sufficiently regular solutions of $i\partial_tu=A_3u$, the norm $\|u(t)\|_{H_3}$ and the energy $q_3(u(t),u(t))$ are constant.
\end{corollary}
\begin{proof}
The Reynolds projection commutes with the functional calculus by Theorem~\ref{thm:sheetoperator}. Unitarity preserves the norm, and commutation with $A_3^{1/2}$ preserves the energy on the form domain.
\end{proof}

\subsection{A precise extension problem}

Let $K$ be a compact truncation of the sheeted core with an external boundary $\partial_{\mathrm{ext}}K$. Connect two cut pieces when they share an identified seam arc of positive length, and assume every component of this adjacency graph meets an external Dirichlet boundary portion of positive length. Let $V(K)$ be the space with matching seam traces, and $V_0(K)$ its subspace with zero external trace.

\begin{lemma}[Coercivity on the compact core]\label{lem:poincare}
There is $C>0$ such that $\|u\|_{L^2(K)}^2\le Cq_K(u,u)$ for every $u\in V_0(K)$.
\end{lemma}
\begin{proof}
Otherwise choose $u_j$ with $L^2$ norm one and energy tending to zero. The sequence is bounded in the finite piecewise $H^1$ space. Rellich compactness gives a subsequence converging strongly in $L^2$ and weakly in $H^1$ to $u$. Its gradient vanishes on every piece. The trace constraints pass to the weak limit, so these constants agree across each connected component. The external zero trace makes each constant zero, contradicting the unit $L^2$ norm.
\end{proof}

\begin{theorem}[Harmonic extension with seam matching]\label{thm:harmonic}
For every admissible external trace $b\in\Tr_{\mathrm{ext}}V(K)$, there is a unique energy-minimizing $u_b\in V(K)$ with trace $b$. It is characterized by
\begin{equation}
q_K(u_b,v)=0\qquad(v\in V_0(K)).
\end{equation}
It is harmonic on the interiors, satisfies the conservative seam condition, and depends continuously on $b$ in the quotient trace norm. The extension map is equivariant whenever the boundary datum is.
\end{theorem}
\begin{proof}
Choose a lift $w\in V(K)$ of $b$. Seek $u_b=w+z$ with $z\in V_0(K)$. The equation for $z$ is $q_K(z,v)=-q_K(w,v)$. Lemma~\ref{lem:poincare} makes $q_K$ coercive on $V_0(K)$, so Lax--Milgram gives a unique $z$ and an energy bound. Replacing $w$ by another lift changes it by an element of $V_0(K)$ and produces the same $u_b$. Taking the infimum over lifts gives continuity in the quotient trace norm. Orthogonality shows
$q_K(u_b+v,u_b+v)=q_K(u_b,u_b)+q_K(v,v)$,
proving minimality and uniqueness. Interior and seam equations follow by testing as in Theorem~\ref{thm:sheetoperator}. Equivariance follows by transforming the weak equation and using uniqueness.
\end{proof}

\begin{remark}[What extension adds]
An independent dilation on a sheet can evolve a local profile without using its neighbors. Harmonic extension uses their common traces and balances their fluxes. Its output therefore depends on the gluing relation through the domain of the energy form. This is the analytic content added by the seam.
\end{remark}

The specification of an operator domain is essential in spectral
problems on singular spaces; Cheeger's analysis provides a foundational
treatment of this issue in singular Riemannian geometry
\cite{Cheeger1983}.

\subsection{The matching defect and the seam response}

\begin{proposition}[Local solutions and their seam defect]\label{prop:defect}
Let $u\in V_3$ satisfy $(-\Delta-\lambda)u_\alpha=0$ on each piece, with conormal derivatives interpreted in the weak trace sense. Its failure to solve the global weak eigenvalue equation is the seam functional
\begin{equation}\label{eq:defect}
\mathfrak d_e(u)=\sum_{\alpha\in I(e)}\kappa_\alpha\partial_{\nu_\alpha}u_\alpha|_e.
\end{equation}
Specifically, for test functions supported away from external boundary,
\begin{equation}
q_3(u,v)-\lambda\langle u,v\rangle_{H_3}
=\sum_e\langle\mathfrak d_e(u),\Tr_ev\rangle.
\end{equation}
Thus $u$ is a global weak state precisely when all these seam functionals vanish.
\end{proposition}
\begin{proof}
Apply Green's identity on each piece. The interior terms vanish by the local equation, and the identified test traces allow the conormal terms to be collected edge by edge. Their sum is~\eqref{eq:defect}. Vanishing is sufficient by the weak definition of $A_3$ and necessary by arbitrary local edge tests, interpreted in the compatible trace space.
\end{proof}

This is a precise seam-supported obstruction statement: it concerns the defect of the chosen differential equation. The divisor of a Mellin transform is a different object, treated in Section~\ref{sec:divisors}.

The variational relationship between boundary traces, solution
operators, and Dirichlet-to-Neumann maps can be formulated abstractly
using boundary pairs associated with quadratic forms
\cite{Post2012}.

\begin{proposition}[Dirichlet-to-Neumann matching]\label{prop:DtN}
Cut a compact core into finitely many pieces and fix homogeneous external boundary conditions. Let $\lambda$ lie outside the Dirichlet spectra of the cut pieces. For an admissible seam trace $h$, solve the piecewise eigenvalue equations and let $\cD_\Gamma(\lambda)h$ be the weighted sum of their conormal derivatives. Then the trace map identifies global eigenfunctions at $\lambda$ with
\begin{equation}
\ker\cD_\Gamma(\lambda).
\end{equation}
For real $\lambda$ in this domain, $\cD_\Gamma(\lambda)$ is symmetric under the trace duality pairing.
\end{proposition}
\begin{proof}
The excluded Dirichlet spectra guarantee existence and uniqueness of each piecewise boundary solution by the resolvent of its Dirichlet realization. A global eigenfunction has a common trace and zero seam defect, hence lies over the kernel. Conversely, a kernel trace produces a unique family of solutions, and Proposition~\ref{prop:defect} glues them to a global state. A zero trace gives only the zero piecewise solution, so the correspondence is injective. Symmetry follows by applying the second Green identity to two boundary solutions at the same real $\lambda$ and summing over pieces.
\end{proof}

\section{Boundary fields and transfer to the spectral channels}\label{sec:panels}
\purpose{Specify how the six decorated panels provide input to the operators already defined.}

On a face $F$, consider the auxiliary tensor
\begin{equation}
h_F=du^2-\mu_F(u,w)dw^2,
\qquad \mu_F=(u-\tfrac12)(w-\tfrac12).
\end{equation}
It is positive definite for $\mu_F<0$, Lorentzian for $\mu_F>0$, and degenerate on $C_F$. It is a boundary decoration, separate from the Riemannian propagation metric.

\begin{lemma}[Equivariant chart decoration]\label{lem:signature}
The cross $C_F$ is preserved by the face dihedral group. The sign of the displayed $\mu_F$ is reversed by a quarter-turn. Consequently the full equivariant decoration is the family of transported tensors with their charts, rather than an invariant scalar sign on the face.
\end{lemma}
\begin{proof}
A quarter-turn may be written $(u,w)\mapsto(1-w,u)$. It exchanges the two central lines and sends $(u-1/2)(w-1/2)$ to its negative. Pullback transports a tensor intrinsically; retaining the orbit of the decorated chart records this action without declaring its coordinate formula fixed.
\end{proof}

For $\sigma>0$, define $\phi_\sigma(q)=\exp(\sigma/\log q)$ on $(0,1)$. Its values lie in $(0,1)$, with limits $1$ at $0$ and $0$ at $1$. The four channels on each face use $q=u,1-u,w,1-w$. A corner $(a,b)\in\{0,1\}^2$ selects
\begin{equation}
\Phi_F^{a,b}(u,w)=\phi_{\sigma_u}(q_a(u))\phi_{\sigma_w}(q_b(w)),
\qquad q_0(t)=t,\quad q_1(t)=1-t.
\end{equation}
These are bounded boundary inputs, smooth in the open face. Parameters and channels are transported along with the charts.

\begin{lemma}[Reciprocal-logarithmic integral]\label{lem:reciprocallog}
For $a,\sigma>0$,
\begin{equation}
\int_0^\infty e^{-\sigma/t-at}\dd t
=2\sqrt{\frac{\sigma}{a}}K_1(2\sqrt{a\sigma}).
\end{equation}
More generally, for $\nu\in\C$,
\begin{equation}\label{eq:Besselintegral}
\int_0^\infty t^{\nu-1}e^{-at-\sigma/t}\dd t
=2(\sigma/a)^{\nu/2}K_\nu(2\sqrt{a\sigma}).
\end{equation}
\end{lemma}
\begin{proof}
Substitute $t=\sqrt{\sigma/a}\,e^u$. The exponent becomes $-2\sqrt{a\sigma}\cosh u$, and the power and measure give $(\sigma/a)^{\nu/2}e^{\nu u}\dd u$. Split the integral over $\R$ at zero and combine the two halves into $2\cosh(\nu u)$. This is the integral representation $K_\nu(z)=\int_0^\infty e^{-z\cosh u}\cosh(\nu u)\dd u$ for $z>0$. Both endpoints are exponentially controlled, justifying the substitution for every fixed $\nu$; see also~\cite{DLMFK}.
\end{proof}

The substitution $t=-\log q$ turns the boundary channel into $e^{-\sigma/t}$. Its $K_1$ integral and the $K_{ir}$ modes below belong to the same Bessel family, with their orders determined by different operations.

\begin{definition}[Boundary transfer]
A transfer to the eight cusp collars is a linear map
\begin{equation}
\cR:\cB_{\partial X}\longrightarrow\bigoplus_{v\in V}C^\infty(S^1),
\end{equation}
where $\cB_{\partial X}$ is a specified space of decorated panel fields. A concrete choice integrates corner fields against smooth kernels on the face variables and the output circle. Mean-zero projection $Q_0b=b-\int_0^1b(x)\dd x$ selects the nonconstant channel.
\end{definition}

\begin{proposition}[Equivariant transfer by averaging]\label{prop:transfer}
For representations $U^{\mathrm{in}}$ and $U^{\mathrm{out}}$ of $\G$, the average
\begin{equation}
\cR_{\G}=\frac1{48}\sum_{g\in\G}
U_g^{\mathrm{out}}\cR_0(U_g^{\mathrm{in}})^{-1}
\end{equation}
intertwines the two actions. If the output action preserves the circle measure, mean-zero projection also intertwines it.
\end{proposition}
\begin{proof}
Conjugating the sum by the actions of an element $h$ merely reindexes it by $g\mapsto hg$. Measure preservation makes the integral defining the constant mode invariant.
\end{proof}

A transfer is part of the model. Its role is to specify how a field on a two-dimensional panel supplies data on a one-dimensional trace circle; incidence alone supplies the labels but does not determine this sampling operation.

\section{Cusp Fourier analysis and propagation}\label{sec:cusp}
\purpose{Solve the local spectral equation at each of the eight marked ends and fix its boundary normalization.}

The cusp calculations belong to the spectral analysis of hyperbolic
surfaces, whose geometric and scattering aspects are developed in
\cite{Borthwick2016}.

Use the cusp metric~\eqref{eq:cuspmetric} and the positive Laplacian
\begin{equation}
A_C=-y^2(\partial_x^2+\partial_y^2).
\end{equation}
A collar is intrinsically the quotient of a horoball by $z\mapsto z+1$. The quotient creates the closed circular direction.

\begin{theorem}[Fourier--Bessel cusp modes]\label{thm:cusp}
Let $r\in\R$ and suppose $A_C\Psi=(1/4+r^2)\Psi$. For each nonzero Fourier index $n$, the solution that decays as $y\to\infty$ is a scalar multiple of
\begin{equation}\label{eq:cuspwave}
W_{n,r}(x,y)=\sqrt yK_{ir}(2\pi|n|y)e^{2\pi inx}.
\end{equation}
The zero mode is a linear combination of $y^{1/2+ir}$ and $y^{1/2-ir}$ for $r\ne0$, and of $y^{1/2}$ and $y^{1/2}\log y$ for $r=0$.
\end{theorem}
\begin{proof}
Write $\Psi=\sum_na_n(y)e^{2\pi inx}$. For $n\ne0$ the radial equation is
\begin{equation}
a_n''+\left(\frac{1/4+r^2}{y^2}-4\pi^2n^2\right)a_n=0.
\end{equation}
Set $a_n(y)=\sqrt y\,b(2\pi|n|y)$. Direct differentiation gives
$z^2b''+zb'-(z^2+(ir)^2)b=0$.
The decaying solution of this modified Bessel equation is $K_{ir}(z)$; the independent solution has exponential growth. For $n=0$, substituting $a_0=y^\eta$ gives $\eta(1-\eta)=1/4+r^2$. The repeated-root case yields the logarithmic second solution.
\end{proof}

\begin{proposition}[Normalized collar extension]\label{prop:collar}
Fix $r\in\R$ and choose $Y$ sufficiently large that $K_{ir}(2\pi nY)\ne0$ for all $n\ge1$. For a smooth mean-zero trace $b(x)=\sum_{n\ne0}b(n)e^{2\pi inx}$, the decaying collar solution is
\begin{equation}\label{eq:collarsolution}
\Psi_b(x,y)=\sum_{n\ne0}b(n)
\frac{\sqrt yK_{ir}(2\pi|n|y)}{\sqrt YK_{ir}(2\pi|n|Y)}e^{2\pi inx},
\qquad y\ge Y.
\end{equation}
It has trace $b$ at $Y$. On each region $y\ge Y+\varepsilon$, the series and all derivatives converge normally.
\end{proposition}
\begin{proof}
For fixed $r$, the asymptotic $K_{ir}(z)=\sqrt{\pi/(2z)}e^{-z}(1+O(z^{-1}))$ gives a positive nonzero value for all sufficiently large positive $z$. This permits a common $Y$. The ratio in~\eqref{eq:collarsolution} has exponential factor $e^{-2\pi|n|(y-Y)}$; its differentiated forms have at most additional polynomial factors in $|n|$ on compact $y$-ranges. Smoothness of $b$ gives rapid Fourier decay, and the exponential factor gives normal convergence away from the initial collar. At $y=Y$ the ratio is one. Modewise uniqueness of the decaying solution proves uniqueness among zero-mode-free decaying solutions with that trace.
\end{proof}

\begin{remark}[Trace coefficients and wave coefficients]
In the Fourier--Bessel expansion, the wave coefficients are
\begin{equation}
A(n)=\frac{b(n)}{\sqrt YK_{ir}(2\pi|n|Y)}.
\end{equation}
A smooth collar trace can therefore lead to exponentially growing $A(n)$. A full Mellin integral down to $y=0$ requires an additional global coefficient condition. The local formula~\eqref{eq:collarsolution} is valid on its stated collar without that condition.
\end{remark}

\subsection{Scale, principal series, and dilation}

The cusp eigenvalue can be written $\eta(1-\eta)$ with $\eta=1/2+ir$. Reflection $\eta\mapsto1-\eta$ preserves this eigenvalue. This line labels the oscillatory principal channel. The complete self-adjoint sheet operator may also have other discrete spectral values; the notation describes a channel rather than an assumption about its entire spectrum.

On $L^2(\R_+,dt/t)$, dilation $(D_af)(t)=f(at)$ is unitary. With $x=\log t$, it becomes translation on $L^2(\R,dx)$. Thus the Mellin transform is a Fourier transform in logarithmic scale after the appropriate power normalization. This is why periodic geometry supplies a discrete Fourier label while noncompact scale supplies a Mellin variable.

\section{The ambient hyperbolic realization}\label{sec:bulk}
\purpose{Explain the boundary-to-bulk picture and specify how it couples to the intrinsic sheet equations.}

Recenter the cube as $[-1,1]^3$. Radial projection $q\mapsto q/\|q\|$ identifies its boundary with the ideal sphere of the Poincar\'e ball
\begin{equation}
\Hh^3=\{X\in\R^3:|X|<1\},\qquad
 g_{\Hh^3}=\frac{4\,dX^2}{(1-|X|^2)^2}.
\end{equation}
The eight corners determine the directions
$\xi_\epsilon=(\epsilon_1,\epsilon_2,\epsilon_3)/\sqrt3$,
with $\epsilon_j=\pm1$. Place six seam vertices at $\pm\rho e_j$, for $0<\rho<1$, and join nonopposite pairs by hyperbolic geodesic segments. The eight octahedral triangles and the eight ideal directions carry matching sign labels. This realizes the seam inside the ball and gives a spatial arrangement for the marked ends.

The intrinsic cusp circles come from parabolic quotients. Their developing charts may be drawn toward these ideal directions, while the global quotient metrics and seam attachments remain part of the sheet datum.

\begin{theorem}[Hyperbolic Poisson transform]\label{thm:bulk}
For $f\in L^1(S^2)$ and $\beta\in\C$, define
\begin{equation}\label{eq:bulkPoisson}
(\cH_\beta f)(X)=\int_{S^2}
\left(\frac{1-|X|^2}{|X-\xi|^2}\right)^\beta f(\xi)\dd\omega(\xi).
\end{equation}
Then $\cH_\beta f$ is smooth in the ball and
\begin{equation}
A_{\Hh^3}\cH_\beta f=\beta(2-\beta)\cH_\beta f.
\end{equation}
The transform is equivariant for the orthogonal action preserving the ball, including $\Pi_3$.
\end{theorem}
\begin{proof}
Let $B_\xi(X)=\log((1-|X|^2)/|X-\xi|^2)$. In upper half-space coordinates with $\xi=\infty$, this is $\log y$ up to an additive constant. The hyperbolic identities are $|\nabla B_\xi|^2=1$ and $\Delta B_\xi=-2$. Therefore
$-\Delta e^{\beta B_\xi}=\beta(2-\beta)e^{\beta B_\xi}$.
Isometries give the same identity for every ideal point. On a compact subset of the ball the kernel and all its derivatives are uniformly bounded in $\xi$, so differentiation under the $L^1$ integral is justified. Orthogonal transformations preserve its distance formula and the spherical measure, proving equivariance.
\end{proof}

For asymptotically hyperbolic realizations, the relevant analytic
background includes meromorphic continuation of the resolvent
\cite{MazzeoMelrose1987} and the relation between scattering operators
and conformal boundary geometry \cite{GrahamZworski2003}.
The precise continuation properties depend on the asymptotic metric;
the role of its evenness is analyzed in \cite{Guillarmou2005}.

\begin{lemma}[The trace equation on a geodesic plane]\label{lem:bulktrace}
Let $N\cong\Hh^2$ be a totally geodesic plane in $\Hh^3$, let $r$ be signed normal distance, and let $U$ be a smooth bulk eigenfunction with $A_{\Hh^3}U=\lambda U$. Its trace $u=U|_N$ satisfies
\begin{equation}
A_Nu=\lambda u+\partial_r^2U|_{r=0}.
\end{equation}
\end{lemma}
\begin{proof}
The normal-coordinate metric is $dr^2+\cosh^2r\,g_N$, so
$\Delta_{\Hh^3}=\partial_r^2+2\tanh r\,\partial_r+\cosh^{-2}r\,\Delta_N$.
Evaluate at $r=0$ and use the positive-Laplacian convention.
\end{proof}

The bulk trace can therefore supply boundary data or a forcing term for a sheet problem. An intertwining transfer is a stronger choice: it prescribes a relation between bulk and leafwise eigenspaces. In particular, the bulk principal parameter $\beta=1+i\varrho$ and the cusp parameter $\eta=1/2+ir$ remain separate until a coupling relates them.

The Poisson kernel is also the scalar boundary-to-bulk kernel of Euclidean anti-de Sitter geometry~\cite{Witten}. This supplies a precise holographic aspect of the construction: boundary fields determine interior fields through a transform. The additional seam, its mixed-type chart decoration, and its transmission conditions specify further interior structure.

\section{Automorphic specialization and the Maass Mellin transform}\label{sec:maassmellin}
\purpose{Identify the global input that turns local cusp modes into completed automorphic states with a functional equation.}

We now impose the automorphic hypotheses needed for the classical
Maass-form construction.
Fourier expansions at cusps and the associated spectral methods
are treated systematically in \cite{Iwaniec2002}.

The translations $z\mapsto z+a$, dilations $z\mapsto e^t z$, and inversion $S:z\mapsto-1/z$ lie in $\PSL(2,\R)$. The unit translation $T:z\mapsto z+1$ creates the cusp period, and $S,T$ generate $\PSL(2,\Z)$. A global automorphic model includes an actual discrete group and its transformation law. The periodicity of a collar supplies its parabolic identification; the remaining group relations supply additional global information.

\begin{lemma}[Mellin integral of a Bessel mode]\label{lem:BesselMellin}
For $\nu\in\C$, $c>0$, and $\Re s>|\Re\nu|$,
\begin{equation}\label{eq:BesselMellin}
\int_0^\infty K_\nu(cy)y^{s-1}\dd y
=2^{s-2}c^{-s}\Gamma\!\left(\frac{s+\nu}{2}\right)
\Gamma\!\left(\frac{s-\nu}{2}\right).
\end{equation}
\end{lemma}
\begin{proof}
Use the integral representation
$K_\nu(z)=\tfrac12(z/2)^\nu\int_0^\infty a^{-\nu-1}e^{-a-z^2/(4a)}\dd a$,
which follows by substitution from~\eqref{eq:Besselintegral}. Integrating first in $z$ gives
\begin{equation}
\int_0^\infty z^{s+\nu-1}e^{-z^2/(4a)}\dd z
=2^{s+\nu-1}a^{(s+\nu)/2}\Gamma((s+\nu)/2).
\end{equation}
The remaining $a$ integral is $\Gamma((s-\nu)/2)$. The condition on $\Re s$ makes both absolute integrals finite, justifying Fubini. Finally substitute $z=cy$. This is the standard identity recorded in~\cite{DLMFMellin}.
\end{proof}

For $\nu=ir$ and $c=2\pi n$, the right side is
$\tfrac14\pi^{-s}n^{-s}\Gamma((s+ir)/2)\Gamma((s-ir)/2)$.
Thus a coefficient sequence of polynomial growth has, in a sufficiently far right half-plane, the completed positive-frequency state
\begin{equation}\label{eq:genericMaasscompletion}
4\int_0^\infty\sum_{n\ge1}A(n)K_{ir}(2\pi ny)y^{s-1}\dd y
=\pi^{-s}\Gamma\!\left(\frac{s+ir}{2}\right)
\Gamma\!\left(\frac{s-ir}{2}\right)\sum_{n\ge1}\frac{A(n)}{n^s}.
\end{equation}
Absolute convergence follows by the same Bessel integral at a real exponent exceeding the coefficient growth. A local collar instead gives the integral from $Y$ to infinity, an incomplete Mellin kernel.

\begin{theorem}[Full-level Maass functional equation]\label{thm:maassFE}
Let $u$ be a weight-zero Maass cusp form for $\PSL(2,\Z)$ with eigenvalue $1/4+r^2$, $r\in\R$, and parity $\delta\in\{0,1\}$:
\begin{equation}
u(-x+iy)=(-1)^\delta u(x+iy),\qquad
u(z)=\sqrt y\sum_{n\ne0}A(n)K_{ir}(2\pi|n|y)e^{2\pi inx}.
\end{equation}
Assume the usual polynomial growth of its Fourier coefficients. Then
\begin{equation}
\Lambda_u(s)=\pi^{-(s+\delta)}
\Gamma\!\left(\frac{s+\delta+ir}{2}\right)
\Gamma\!\left(\frac{s+\delta-ir}{2}\right)
\sum_{n\ge1}A(n)n^{-s}
\end{equation}
extends to an entire function and satisfies
\begin{equation}\label{eq:MaassFE}
\Lambda_u(s)=(-1)^\delta\Lambda_u(1-s).
\end{equation}
\end{theorem}
\begin{proof}
For even parity, set $I_0(s)=\int_0^\infty u(iy)y^{s-1/2}\dd y/y$. Parity gives $A(-n)=A(n)$, and Lemma~\ref{lem:BesselMellin} yields $I_0(s)=\tfrac12\Lambda_u(s)$ initially in a right half-plane. Invariance under inversion gives $u(i/y)=u(iy)$. Substituting $y\mapsto1/y$ therefore gives $I_0(s)=I_0(1-s)$.

For odd parity let $v(y)=\partial_xu(x+iy)|_{x=0}$. Then
$v(y)=4\pi i\sqrt y\sum_{n\ge1}nA(n)K_{ir}(2\pi ny)$.
Define $I_1(s)=\int_0^\infty v(y)y^{s+1/2}\dd y/y$. The Bessel integral gives $I_1(s)=i\pi\Lambda_u(s)$. Differentiating $u(z)=u(-1/z)$ at $z=iy$ gives $v(y)=-y^{-2}v(1/y)$, so substitution yields $I_1(s)=-I_1(1-s)$.

At infinity, cusp decay controls both $u$ and its derivative exponentially. The inversion relations transfer this to decay faster than any power at zero. The two period integrals and all their $s$-derivatives therefore converge locally uniformly for every $s\in\C$, proving entire continuation and the stated functional equations.
\end{proof}

\begin{remark}[Arithmetic and analytic inputs]
The theorem uses global inversion invariance. Its gamma factors come from the Bessel integral; its reflection comes from the automorphy law. A simultaneous Hecke eigenstate supplies the multiplicative coefficient relations needed for an Euler product. These are separate structural inputs in the automorphic specialization~\cite{MillerSchmid}.
\end{remark}

\begin{remark}[Two reflections]
The transformation $\eta\mapsto1-\eta$ leaves the local eigenvalue $\eta(1-\eta)$ unchanged. The transformation $s\mapsto1-s$ in~\eqref{eq:MaassFE} acts on an independent Mellin variable and follows from inversion of a global period. The shared form of the reflections is expressed through the transform, with their variables distinguished.
\end{remark}

\section{Higher-rank theta states and normalized duality}\label{sec:higher}
\purpose{Extend the transform mechanism while keeping geometric dimension, lattice rank, and spectral normalization distinct.}

The higher seed notation is
\begin{equation}\label{eq:higherseed}
\mathfrak S_n=(\cI_n,\Gamma_n,\Pi_n),\qquad
\Pi_n\cong B_n,\qquad
\cI_n=\left(\bigsqcup_{i=1}^{2^{n-1}}\Sigma^{n-3}X_i\right)\Big/\!\sim_{\Gamma_n},
\quad n\ge3.
\end{equation}
The groups $B_n$ are the finite Coxeter groups of type $B$,
realized here as signed permutation groups \cite{Humphreys1990}. Here $\Sigma$ denotes the specified suspension operation on sheets; the attaching maps are part of the datum. At $n=3$ there are four unsuspended two-dimensional sheets. The formal sheet dimension at level $n$ is $n-1$. The theta rank is the separate parameter
\begin{equation}
k_n=n-2.
\end{equation}
The action of $B_n$ on geometry and its action on a rank-$k_n$ lattice or kernel space must be related by specified equivariant maps. The two ranks have different meanings.

\subsection{Lattice theta states}

\begin{definition}[Normalized lattice pair]
Let $L\subset\R^k$ be a full lattice of covolume $v>0$, and let
$L^*=\{\xi:\langle\xi,\ell\rangle\in\Z\text{ for all }\ell\in L\}$.
Set
\begin{equation}
\theta_L(t)=\sum_{\ell\in L}e^{-\pi t|\ell|^2},\qquad
\Theta_L(t)=v^{1/2}\theta_L(t).
\end{equation}
The paired constants at infinity are $a_L=v^{1/2}$ and $a_{L^*}=v^{-1/2}$.
\end{definition}

\begin{proposition}[Lattice reciprocity]\label{prop:latticereciprocity}
For every $t>0$,
\begin{equation}\label{eq:latticereciprocity}
\theta_L(t)=v^{-1}t^{-k/2}\theta_{L^*}(1/t),
\qquad
\Theta_L(t)=t^{-k/2}\Theta_{L^*}(1/t).
\end{equation}
\end{proposition}
\begin{proof}
Write $L=B\Z^k$ with $|\det B|=v$. Periodizing a Schwartz function on $\R^k/L$ as in Lemma~\ref{lem:poisson} gives
$\sum_{\ell\in L}f(\ell)=v^{-1}\sum_{\ell^*\in L^*}\widehat f(\ell^*)$.
The $k$-dimensional Gaussian transform is the product of the one-dimensional transforms, equal to $t^{-k/2}e^{-\pi|\xi|^2/t}$. This proves the first formula. Since $L^*$ has covolume $v^{-1}$, multiplying by $v^{1/2}$ proves the second.
\end{proof}

Our lattice and dual-lattice conventions follow the standard
Euclidean framework for theta series
\cite{ConwaySloane1999}.

\begin{definition}[Admissible reciprocal theta pair]\label{def:thetapair}
A reciprocal theta pair of weight $k>0$ consists of smooth kernels $\theta_\pm:(0,\infty)\to\C$ and constants $a_\pm$ such that
\begin{equation}\label{eq:pairedreciprocity}
\theta_+(t)=t^{-k/2}\theta_-(1/t),
\end{equation}
and $\theta_\pm(t)-a_\pm$ decay exponentially as $t\to\infty$. Families of such pairs are required to have locally uniform decay bounds in their other parameters whenever parameterwise analytic continuation is used.
\end{definition}

Normalized dual lattices provide examples. The definition also allows kernels built by other equivariant spectral constructions, as long as their constants and reciprocity are specified.

\begin{theorem}[Mellin reflection for a reciprocal pair]\label{thm:pairedMellin}
For $\Re s>k$, define
\begin{equation}
\Lambda_\pm(s)=\frac12\int_0^\infty(\theta_\pm(t)-a_\pm)t^{s/2}\frac{\dd t}{t}.
\end{equation}
These functions extend meromorphically to $\C$ by
\begin{align}\label{eq:pairedMellin}
\Lambda_\pm(s)={}&-\frac{a_\pm}{s}+\frac{a_\mp}{s-k}\\
&+\frac12\int_1^\infty\left[(\theta_\pm(t)-a_\pm)t^{s/2}
+(\theta_\mp(t)-a_\mp)t^{(k-s)/2}\right]\frac{\dd t}{t}.\nonumber
\end{align}
They satisfy $\Lambda_+(s)=\Lambda_-(k-s)$. Their only possible poles are at $0,k$, with residues $-a_\pm,a_\mp$, respectively. The regularized functions
$\Xi_\pm(s)=s(s-k)\Lambda_\pm(s)$ are entire and obey the same paired reflection.
\end{theorem}
\begin{proof}
At zero, reciprocity gives
$\theta_\pm(t)=t^{-k/2}(a_\mp+O(e^{-c/t}))$,
so the initial integral converges absolutely for $\Re s>k$. Split it at $1$ and use reciprocity in the small-$t$ part. The terms involving the constants integrate to $a_\mp/(s-k)-a_\pm/s$, while the other term changes to the second integral over $[1,\infty)$. Exponential decay makes that integral entire by locally uniform domination, including all derivatives in $s$. Replacing $s$ by $k-s$ and exchanging the signs leaves the expression unchanged. The pole statements and entire regularization follow directly.
\end{proof}

\begin{corollary}[A common reflection coordinate]\label{cor:normalizedw}
At weight $k$, put $s=kw$ and $\widetilde\Lambda_\pm(w)=\Lambda_\pm(kw)$. Then
\begin{equation}
\widetilde\Lambda_+(w)=\widetilde\Lambda_-(1-w).
\end{equation}
The residues at $w=0,1$ are $-a_\pm/k,a_\mp/k$.
\end{corollary}
\begin{proof}
Substitute $s=kw$ into Theorem~\ref{thm:pairedMellin}. A linear change of variable divides each residue by $k$.
\end{proof}

The common reflection coordinate is useful in the tower, but it is a normalization choice. For an isotropic rank-$k$ theta kernel, the natural unscaled reflection is $s\mapsto k-s$. A multivariable product model instead has one independent reflection for each scale coordinate. These two Mellin constructions retain different information.

\subsection{Finite symmetry and theta data}

\begin{proposition}[Equivariant averaging commutes with Mellin transformation]\label{prop:avgMellin}
Let a finite group act linearly on a family of reciprocal theta pairs, preserving their weight and the stated analytic bounds. Transport the constants with the kernels. Then Reynolds averaging preserves reciprocity and commutes with the regularized Mellin construction~\eqref{eq:pairedMellin}.
\end{proposition}
\begin{proof}
The reciprocity equation is linear in the two kernels. Formula~\eqref{eq:pairedMellin}, including its constant terms, is also linear. A finite average can be interchanged with its convergent integrals and rational terms.
\end{proof}

This is the precise operation used when local theta channels are first made compatible with geometry and then projected to an invariant sector. The existence and choice of those local channels belong to the spectral enrichment; the commutation statement says that their finite symmetry reduction can be performed before or after Mellin transformation.

\section{Changing rank: marginalization, residues, and duality}\label{sec:ranktransition}
\purpose{Give a concrete rank-two to rank-one calculation and expose the compatibility that a geometric bonding map must induce on spectral data.}

\subsection{What Fourier duality does to a deleted coordinate}

For $f\in\Sch(\R^{k+1})$, define integration and restriction in the final coordinate by
\begin{equation}
(\rho f)(x')=\int_\R f(x',u)\dd u,
\qquad (\pi f)(x')=f(x',0).
\end{equation}
Let $\mathcal F_k$ be the Fourier transform on $\R^k$ with the convention of Lemma~\ref{lem:poisson}.

\begin{lemma}[The Fourier exchange law]\label{lem:marginal}
On Schwartz spaces,
\begin{equation}\label{eq:Fourierexchange}
\pi\mathcal F_{k+1}=\mathcal F_k\rho,
\qquad
\rho\mathcal F_{k+1}=\mathcal F_k\pi.
\end{equation}
\end{lemma}
\begin{proof}
Evaluating $\widehat f(\xi',\eta)$ at $\eta=0$ gives the Fourier transform in $x'$ of $\int f(x',u)\dd u$, by Fubini. For the second identity, apply one-dimensional Fourier inversion in $u$ at $u=0$, then Fourier transform in $x'$. Schwartz decay justifies all operations.
\end{proof}

Thus an integration map changes under duality into a restriction map. A rank-lowering prescription has to record which of these operations is being used on each side.

\begin{proposition}[A compatible Gaussian lift]\label{prop:Gaussianlift}
Let $g(u)=e^{-\pi u^2}$ and $Ef=f\otimes g$. Then
\begin{equation}
\rho E=\pi E=\mathrm{id},\qquad
\mathcal F_{k+1}E=E\mathcal F_k.
\end{equation}
On the image of $E$, integration is a Fourier-compatible inverse to the lift.
\end{proposition}
\begin{proof}
The Gaussian satisfies $g(0)=\int_\R g=1$ and $\widehat g=g$. Tensor factorization of the Fourier transform proves all three identities.
\end{proof}

For periodic heat kernels, the corresponding geometric operation is an angular constant term. Define
\begin{equation}
K_k(\mathbf t,\mathbf x)=\sum_{\mathbf m\in\Z^k}
\exp\left(-\pi\sum_{j=1}^kt_jm_j^2\right)e^{2\pi i\mathbf m\cdot\mathbf x}.
\end{equation}
For $t_j>0$, absolute convergence gives
\begin{equation}\label{eq:heatconstantterm}
\int_0^1K_{k+1}(\mathbf t,\mathbf x)\dd x_{k+1}
=K_k((t_1,\ldots,t_k),(x_1,\ldots,x_k)).
\end{equation}
The extra circular Fourier index is projected to zero.

\subsection{An exact split rank-two Mellin calculation}

Let $h(t)=(\vartheta(t)-1)/2=\sum_{m\ge1}e^{-\pi m^2t}$, where $\vartheta=\theta_\Gamma$. In the split rank-two model use two independent heat parameters and the kernel $h(t_1)h(t_2)$. Its relation to the full rank-two theta kernel is
\begin{equation}
4h(t_1)h(t_2)=\vartheta(t_1)\vartheta(t_2)-\vartheta(t_1)-\vartheta(t_2)+1.
\end{equation}
This removes the two coordinate-axis constant sectors before the double Mellin transform.

\begin{theorem}[Split rank-two Mellin compatibility]\label{thm:ranktwo}
For $\Re s_1,\Re s_2>1$,
\begin{equation}\label{eq:ranktwointegral}
\mathcal L_2(s_1,s_2):=\int_0^\infty\!\int_0^\infty
h(t_1)h(t_2)t_1^{s_1/2}t_2^{s_2/2}
\frac{\dd t_1}{t_1}\frac{\dd t_2}{t_2}
=\Lam(s_1)\Lam(s_2).
\end{equation}
After meromorphic continuation,
\begin{equation}\label{eq:ranktworesidue}
\Res_{s_2=1}\mathcal L_2(s_1,s_2)=\Lam(s_1),\qquad
-\Res_{s_2=0}\mathcal L_2(s_1,s_2)=\Lam(s_1).
\end{equation}
With $R_2F(s_1,s_2)=F(1-s_1,1-s_2)$, $R_1f(s)=f(1-s)$,
$U_+=\Res_{s_2=1}$, and $U_-=-\Res_{s_2=0}$, the dual transition law is
\begin{equation}\label{eq:residueintertwine}
U_-R_2=R_1U_+.
\end{equation}
\end{theorem}
\begin{proof}
Absolute convergence in the product half-plane permits Fubini, so the double integral is the product of the two integrals in Theorem~\ref{thm:seamzeta}. Their meromorphic continuations give~\eqref{eq:ranktwointegral} globally. The residues of $\Lam$ at $1$ and $0$ are $1$ and $-1$, giving~\eqref{eq:ranktworesidue}.

For the last identity, expand $F(1-s_1,z)$ near $z=1$ as $a_{-1}(1-s_1)/(z-1)$ plus the remaining Laurent terms. Substituting $z=1-s_2$ changes this term to $-a_{-1}(1-s_1)/s_2$. Applying $-\Res_{s_2=0}$ yields $a_{-1}(1-s_1)$, which is $R_1U_+F$. The argument applies coefficientwise wherever the relevant meromorphic Laurent germs are defined, and hence by continuation in the remaining variable.
\end{proof}

\begin{corollary}[A product tower]
The functions $\mathcal L_k(s_1,\ldots,s_k)=\prod_{j=1}^k\Lam(s_j)$ satisfy independent functional equations in each coordinate. Taking the residue at $s_k=1$ lowers the rank by one; the dual prescription is minus the residue at $s_k=0$.
\end{corollary}
\begin{proof}
Apply the rank-one reflection to each factor and repeat Theorem~\ref{thm:ranktwo} in the final coordinate.
\end{proof}

\begin{remark}[The meaning of the model transition]
At the level $\mathfrak S_4\to\mathfrak S_3$, theta rank drops from two to one. The theorem gives a complete analytic model for that drop when the added spectral direction is the split Gaussian factor. To realize it for a geometric bonding morphism, its action on the chosen heat or theta states must be this specified reduction. The theorem identifies the exact Mellin transition that such a realization must induce.
\end{remark}

\begin{remark}[Residues retain more than a boundary value]
A residue is the coefficient of a polar germ. It is not evaluation on the line $\Re s=1$. The distinction permits a precise version of gluing near singular spectral strata: a transition acts on meromorphic states by a stated residue map, and duality exchanges its two polar prescriptions. The isotropic Epstein transform obtained by setting $t_1=t_2$ is a different transform and is not the product in~\eqref{eq:ranktwointegral}.
\end{remark}

\section{Enriched bonding morphisms and the doubled base}\label{sec:base}
\purpose{Define the geometric tower without hiding the maps that make its inverse limits and duality meaningful.}

\subsection{Descending sheet maps through the seam}

\begin{definition}[Enriched bonding datum]\label{def:bonding}
A bonding datum $p_n:\mathfrak S_{n+1}\to\mathfrak S_n$ includes the maps of sheet and seam strata, their compatibility with the attaching relations, the specified relation between the symmetry actions, and the induced maps on the chosen analytic data. In a stable formulation, suspension and desuspension are performed in the declared stable category; their geometric realizations, when used, are separately specified.
\end{definition}

\begin{lemma}[Quotient descent criterion]\label{lem:quotientdescent}
Let $q_i:Y_i\to Y_i/\!\sim_i$ be quotient maps and let $f:Y_1\to Y_2$ be continuous. If $y\sim_1y'$ implies $f(y)\sim_2f(y')$, there is a unique continuous map $\overline f$ satisfying $\overline f q_1=q_2f$. If $f$ is equivariant for a specified group homomorphism and the relations are preserved by the actions, $\overline f$ is equivariant as well.
\end{lemma}
\begin{proof}
Define $\overline f([y])=[f(y)]$. The implication makes it well-defined, and surjectivity of $q_1$ gives uniqueness. The quotient topology makes it continuous because its composite with $q_1$ is continuous. Equivariance follows by evaluating the identity on representatives.
\end{proof}

A longitudinal desuspension of individual sheets supplies part of the input $f$. Lemma~\ref{lem:quotientdescent} states the remaining condition needed to descend it through the seam.

\begin{proposition}[Marked-coordinate symmetry reduction]\label{prop:markedcoordinate}
The stabilizer in $B_{n+1}$ of the last coordinate axis is isomorphic to $B_n\times C_2$ and has a natural homomorphism onto $B_n$. Forgetting the last coordinate is not a map from all signed permutation matrices in $B_{n+1}$ to signed permutation matrices in $B_n$.
\end{proposition}
\begin{proof}
An axis-preserving signed permutation has a signed $n\times n$ permutation block and an independent sign in the last coordinate. Projection onto the first block is a homomorphism. In contrast, a permutation exchanging the first and last axes has, after deleting the last row and column, a zero first row and column; the remaining matrix is not a signed permutation matrix.
\end{proof}

One can therefore use a marked-axis stabilizer when defining a projection, or retain all marked directions and their transport as an equivariant correspondence. An inverse system of the groups $\Pi_n$ themselves requires actual group homomorphisms. The objectwise notation~\eqref{eq:higherseed} does not supply those homomorphisms on its own.

\subsection{Two compatible towers}

For the inverse-limit construction, work with specified topological realizations of enriched seeds and continuous bonding maps. Suppose there are positive and negative towers $\mathfrak S_n^\pm$ and homeomorphisms
\begin{equation}
d_n:\mathfrak S_n^+\longrightarrow\mathfrak S_n^-,\qquad
p_n^-d_{n+1}=d_np_n^+.
\end{equation}
The inverse direction is $d_n^{-1}$. Assume that the common three-seed has compatible embeddings $j_n^\pm:\mathfrak S_3\to\mathfrak S_n^\pm$ with
\begin{equation}\label{eq:seedlifts}
p_n^\pm j_{n+1}^\pm=j_n^\pm,\qquad
 d_nj_n^+=j_n^-.
\end{equation}
The embeddings at level three identify the two copies of the central seed.

When the tower is interpreted in an equivariant homotopy or stable
category, its suspensions and bonding morphisms must be specified
in that category; we use the standard equivariant framework
described in \cite{May1996}. The assembly uses the universal properties of limits and pushouts,
together with compatibility expressed by commutative diagrams;
the categorical conventions are those of \cite{MacLane1998}.

\begin{definition}[Doubled inverse-limit base]\label{def:doubledbase}
Define
\begin{equation}
\widehat{\mathfrak S}_\pm=\varprojlim_{n\ge3}\mathfrak S_n^\pm,
\qquad
\LL=\widehat{\mathfrak S}_-\sqcup_{\mathfrak S_3}\widehat{\mathfrak S}_+,
\end{equation}
using the embeddings $j_\pm=(j_n^\pm)_n$ for the pushout. If the seed symbols are used as objects in a category of enriched structures, these formulas refer to limits and pushouts in the stated realization category.
\end{definition}

\begin{theorem}[Geometric duality of the doubled base]\label{thm:baseduality}
The maps $d_n$ induce a homeomorphism $\widehat d:\widehat{\mathfrak S}_+\to\widehat{\mathfrak S}_-$. It descends to an involution $D:\LL\to\LL$ which exchanges the two towers and fixes the common embedded three-seed pointwise. If the level spaces are nonempty compact Hausdorff spaces, their bonding maps are surjective, and the seed embeddings are closed, then the inverse limits and $\LL$ are nonempty compact Hausdorff spaces.
\end{theorem}
\begin{proof}
A compatible sequence $x=(x_n)$ maps to $(d_nx_n)$, which is compatible by the bonding identity. Applying the inverse maps gives the inverse sequence. Both maps are continuous because their coordinates are continuous in the product topologies. Equation~\eqref{eq:seedlifts} makes $\widehat d$ identify the positive seed embedding with the negative one, so it descends through the pushout. Define $D$ by $\widehat d$ on the positive side and $\widehat d^{-1}$ on the negative side. Then $D^2=\mathrm{id}$ and every identified seed point is fixed.

For the topological statement, an inverse limit is the closed subset of the product cut out by $p_nx_{n+1}=x_n$. Compactness follows from Tychonoff. Surjectivity gives the finite intersection property for these constraints, hence nonemptiness. The seed lifts are embeddings because their level-three coordinate is an embedding. Gluing two compact Hausdorff spaces along homeomorphic closed subspaces gives a compact Hausdorff quotient: the equivalence relation is the diagonal together with the two closed graphs of the identification, so it is closed. This proves the assertions.
\end{proof}

\begin{remark}[What is fixed at the center]
The involution fixes the geometric three-seed in the doubled base. Fourier duality can still act nontrivially on spectral data over that fixed point. A self-dual theta state is an additional fiberwise condition, exemplified by the standard Gaussian and by $\theta_\Gamma$.
\end{remark}

\section{Compatible states and the universal Mellin construction}\label{sec:assembly}
\purpose{State the exact hypotheses under which the local transform mechanisms become one compatible object over the doubled geometric base.}

\subsection{Families over a tower}

Let $\Theta_n^\pm\to\mathfrak S_n^\pm$ be families of admissible theta data. A fiber point records the kernel, its dual, its constants, and any transfer or normalization data. Let
$T_n^\pm:\Theta_{n+1}^\pm\to\Theta_n^\pm$ cover $p_n^\pm$. The compatible fiber over $x=(x_n)\in\widehat{\mathfrak S}_\pm$ is
\begin{equation}
\widehat\Theta_x^\pm=\{(\theta_n):\theta_n\in(\Theta_n^\pm)_{x_n},\quad
T_n^\pm\theta_{n+1}=\theta_n\}.
\end{equation}
Existence of a state in this fiber means that such a compatible sequence has been chosen or constructed.

At level $n$, let $W_n$ be its complex spectral parameter space and let $R_n$ be a specified reflection on meromorphic functions. In a single normalized variable, $R_nF(w)=F(1-w)$. In the split model of Section~\ref{sec:ranktransition}, $W_n=\C^{k_n}$ and $R_n$ reflects every coordinate. Write $\Mer(W_n)$ for the space of meromorphic states and
$\cM_n^\pm:\Theta_n^\pm\to\Mer(W_n)$ for the Mellin construction, including its regularization.

\begin{definition}[Mellin-compatible transition]\label{def:Mellintransition}
A transition $U_n^\pm:\Mer(W_{n+1})\to\Mer(W_n)$ is Mellin-compatible with $T_n^\pm$ if
\begin{equation}\label{eq:Mellincompatibility}
\cM_n^\pm T_n^\pm=U_n^\pm\cM_{n+1}^\pm
\end{equation}
on the stated domain. Parameter changes, constant-term subtractions, and residue prescriptions belong to the maps in this identity.
\end{definition}

This is the operational content of compatibility. It may be established in a convergent domain and continued meromorphically. Theorem~\ref{thm:ranktwo} gives an example where $U_n$ is a residue rather than evaluation at a parameter.

\subsection{Duality and assembly}

Suppose fiberwise isomorphisms $J_n:\Theta_n^+\to\Theta_n^-$ cover $d_n$, with inverse maps in the opposite direction, and satisfy
\begin{equation}\label{eq:dualcompatibility}
T_n^-J_{n+1}=J_nT_n^+,
\qquad \cM_n^-J_n=R_n\cM_n^+,
\qquad U_n^-R_{n+1}=R_nU_n^+.
\end{equation}
For reciprocal theta pairs, the middle identity is Theorem~\ref{thm:pairedMellin} in the chosen spectral coordinates. The remaining identities require the transitions to transport that duality consistently.

\begin{theorem}[Duality and assembly]\label{thm:assembly}
Assume the geometric hypotheses defining $D$ in Theorem~\ref{thm:baseduality}, Mellin compatibility~\eqref{eq:Mellincompatibility}, and dual compatibility~\eqref{eq:dualcompatibility}. Then:
\begin{enumerate}
\item A compatible theta state $\theta=(\theta_n)$ over $x$ gives a compatible completed state
\begin{equation}\label{eq:universalstate}
\cZ(x,\theta)=\bigl(\cM_n^\pm\theta_n\bigr)_{n\ge3}.
\end{equation}
\item Duality sends $\theta$ to the compatible state $J\theta=(J_n\theta_n)$ over $Dx$, and
\begin{equation}\label{eq:universalFE}
\cZ(Dx,J\theta)=R\cZ(x,\theta),
\end{equation}
where $R$ acts by $R_n$ on each level.
\item If compatible identifications of the state families and their completed states are supplied over the common three-seed, the two constructions descend through those identifications to the doubled base. A self-dual rank-one state over the common seed obeys $\Lambda(s)=\Lambda(1-s)$.
\end{enumerate}
\end{theorem}
\begin{proof}
For the first assertion,
\begin{equation}
U_n^\pm\cM_{n+1}^\pm\theta_{n+1}
=\cM_n^\pm T_n^\pm\theta_{n+1}
=\cM_n^\pm\theta_n.
\end{equation}
For the second, the first identity of~\eqref{eq:dualcompatibility} gives
$T_n^-J_{n+1}\theta_{n+1}=J_n\theta_n$.
The middle identity gives its Mellin state as $R_n\cM_n^+\theta_n$, and the last identity makes these reflected states compatible with the negative transitions. The inverse maps prove involutivity.

For the third, the supplied identifications make the maps equal on representatives that are to be identified over the seed. They descend by the quotient universal property, as in Lemma~\ref{lem:quotientdescent}. At a fixed geometric point and self-dual fiber state,~\eqref{eq:universalFE} becomes its scalar functional equation. Rank one gives the normalization $s\mapsto1-s$.
\end{proof}

\begin{remark}[The universal section]
Equation~\eqref{eq:universalstate} is intrinsically a compatible family of meromorphic states. It requires no infinite sum or infinite product of their values. In this inverse-system setting it is a pro-meromorphic section. Realization as an ordinary meromorphic section on a specified parameter space requires descent data for that realization. Keeping the compatible family as the basic object accommodates rank-changing residue maps without treating them as scalar pullbacks.
\end{remark}

\begin{lemma}[Descent to a meromorphic line-bundle section]\label{lem:sectiondescent}
Let a complex parameter space $B$ have an open cover $\{B_i\}$ and nonvanishing holomorphic transition functions $g_{ij}$ with $g_{ij}g_{jk}=g_{ik}$. Meromorphic local states satisfying $F_i=g_{ij}F_j$ on overlaps define a unique meromorphic section of the associated line bundle. If a duality carries the cover, transitions, and local functional equations to one another, the section obeys the induced global functional equation.
\end{lemma}
\begin{proof}
Identify $(x,z)_j$ with $(x,g_{ij}(x)z)_i$. The cocycle makes these identifications consistent. The overlap identity makes the local graphs of $F_i$ agree, giving the section and its uniqueness. Applying duality in the local trivializations gives the global identity by the same overlap condition.
\end{proof}

This lemma is the additional step when $\cZ(\mathbf s,\theta)$ denotes a literal meromorphic section on a geometric or analytic parameter space. Theorem~\ref{thm:assembly} handles compatibility through the tower, and Lemma~\ref{lem:sectiondescent} handles its specified local realization.

\subsection{The three-level organization}

At each level, form the family whose fiber over a theta datum is the space of meromorphic spectral states. The Mellin construction selects its graph $\theta\mapsto(\theta,\cM_n\theta)$. Passing to compatible families gives
\begin{equation}\label{eq:finalfibration}
\cZ\longrightarrow\Theta\longrightarrow\LL.
\end{equation}
The right arrow forgets the theta enrichment and retains the geometry. The left arrow forgets the completed state and retains its theta input. These are family projections; a stronger bundle or categorical fibration structure is included when specified. Mellin transformation selects a completed state rather than serving as a downward forgetful arrow.

\begin{example}[A nonempty reference family]
At rank $k$, take $k$ standard Gaussian theta factors with the coordinate-sector subtractions of Section~\ref{sec:ranktransition}. The Mellin state is $\mathcal L_k=\prod_{j=1}^k\Lam(s_j)$. The maps $U_+=\Res_{s_k=1}$ and $U_-=-\Res_{s_k=0}$ produce compatible completed families. Proposition~\ref{prop:Gaussianlift} supplies compatible Gaussian lifts at the Schwartz level. This is an explicit spectral model; a realization over geometric seeds specifies how these factors and reductions arise from their enrichments.
\end{example}

\section{Divisors, seam defects, and spectral localization}\label{sec:divisors}
\purpose{Explain what the completed states say about zeros and poles, and formulate the identity needed for spectral localization.}

A nonzero meromorphic function has a divisor recording the orders of its zeros and poles. For a meromorphic line-bundle section, these orders are defined in local trivializations and are unchanged by nonvanishing holomorphic transition functions. Divisors are therefore natural outputs of the completed-state layer.

\begin{proposition}[Reflection of divisors]\label{prop:divisor}
Let $F$ be a nonzero meromorphic section and let an involution $\iota$ satisfy $\iota^*F=\varepsilon F$, where $\varepsilon$ is nowhere vanishing and holomorphic. Then
\begin{equation}
\iota^*\divisor(F)=\divisor(F).
\end{equation}
In rank one, zeros and poles occur in configurations invariant under $s\mapsto1-s$.
\end{proposition}
\begin{proof}
A nonvanishing holomorphic factor does not change local order. Pullback by a biholomorphism transports each local order to the corresponding point. Apply these facts to the assumed identity.
\end{proof}

\begin{remark}[Invariant divisor and fixed locus]
An invariant divisor need not be supported on the fixed locus: $(s-a)(s-(1-a))$ is invariant under $s\mapsto1-s$ for every $a$. Locating a divisor on a critical line requires information beyond reflection symmetry.
\end{remark}

The seam defect in Proposition~\ref{prop:defect} tests whether local solutions satisfy a global differential equation. A Mellin divisor records zeros and poles of a transformed state that may already be globally compatible. Theorem~\ref{thm:seamzeta}, for example, starts with a well-defined invariant heat state and produces a meromorphic function with its own divisor. A related connection between geometric decomposition and spectral
invariants is provided by determinant gluing formulas.
The analytic formula of Burghelea, Friedlander, and Kappeler
\cite{BurgheleaFriedlanderKappeler1992} and the combinatorial
construction of Reshetikhin and Vertman \cite{ReshetikhinVertman2015}
identify boundary operators as central ingredients in such formulas.
They provide a natural point of comparison for the seam construction.

The regularized Fredholm determinant $\det_2$ is defined for
Hilbert--Schmidt operators, with convergence and analytic properties
governed by trace-ideal theory \cite{Simon2005}.

\begin{proposition}[A self-adjoint determinant criterion]\label{prop:detcriterion}
Let $K$ be compact, self-adjoint, and Hilbert--Schmidt, with nonzero real eigenvalues $\kappa_j$ counted with multiplicity. Then
\begin{equation}
\det{}_2(I-zK)=\prod_j(1-z\kappa_j)e^{z\kappa_j}
\end{equation}
is entire and all its zeros are real. If a zero-free entire function $E$ satisfies
\begin{equation}\label{eq:detcriterion}
E(z)\,\xi(\tfrac12+iz)=\det{}_2(I-zK),
\end{equation}
then every zero of $\xi(s)$ lies on $\Re s=1/2$.
\end{proposition}
\begin{proof}
Hilbert--Schmidt summability gives $\sum_j\kappa_j^2<\infty$. On compact $z$-sets, all sufficiently late factors have $|z\kappa_j|\le1/2$, and
$\log((1-z\kappa_j)e^{z\kappa_j})=O(|z|^2\kappa_j^2)$.
The tail product converges normally to a nonvanishing analytic function; finite initial factors account for its zeros, exactly the real numbers $1/\kappa_j$. The zero-free $E$ does not alter zeros in~\eqref{eq:detcriterion}. A zero $s_0$ of $\xi$ gives a real zero $z_0=(s_0-1/2)/i$, so $\Re s_0=1/2$.
\end{proof}

This criterion makes the spectral localization objective concrete. An identity of the form~\eqref{eq:detcriterion} lets a self-adjoint operator control a zeta divisor. The heat-trace Mellin identity of Theorem~\ref{thm:seamzeta} produces zeta from positive spectral data and has a different conclusion. The framework records which identity a chosen enrichment supplies.

\section{Local systems, colors, and Fano incidence}\label{sec:enrichment}
\purpose{Record the algebraic enrichments of the seam while distinguishing their multiplication laws from geometric symmetry.}

\subsection{Local systems on the seam}

A local system adds parallel transport. A unitary rank-$d$ local system on the graph is described by a unitary representation of $\pi_1(\Gamma_3)$ up to conjugacy. Its edge energy uses covariant derivatives and the chosen fiber identifications at vertices.

\begin{proposition}[Rank-one holonomy parameters]\label{prop:holonomy}
The isomorphism classes of unitary rank-one local systems on $\Gamma_3$ form
\begin{equation}
\Hom(H_1(\Gamma_3;\Z),U(1))\cong U(1)^7.
\end{equation}
An equivariant local system additionally has lifts of the $\Pi_3$ action to its fibers, preserving parallel transport and satisfying the group law.
\end{proposition}
\begin{proof}
A spanning tree has five edges, leaving seven edges outside it. Contracting the tree gives a wedge of seven circles, so the fundamental group is free on seven generators and its abelianization is $\Z^7$. Rank-one unitary holonomy assigns one unit complex number to each generator. Equivariance includes actual fiber maps covering the group action and obeying composition; invariance of an isomorphism class alone does not specify those maps.
\end{proof}

\subsection{The geometric partition \texorpdfstring{$3+4$}{3+4}}

Three opposite-face classes label seam colors, and four body diagonals label sheets. To equip these seven labels with Fano incidence, let $A=\mathbb F_2^2$ be the four-point affine plane and $D=\mathbb F_2^2\setminus\{0\}$ its three directions. Add one point $d_\infty$ for each direction. The lines are
\begin{equation}\label{eq:Fanolines}
\{d_\infty:d\in D\},\qquad
\{a,a+d,d_\infty\}\quad(d\ne0,\ a\in A/\langle d\rangle).
\end{equation}

\begin{theorem}[The $3+4$ incidence correspondence]\label{thm:Fano}
The incidence system~\eqref{eq:Fanolines} is the Fano plane. Its automorphisms preserving the line at infinity form $\operatorname{AGL}(2,2)\cong S_4$. The action of $\Pi_3$ on the four body diagonals factors through $\Pi_3/\{\pm I\}\cong S_4$, and its action on the three opposite-face classes agrees with the action on the three partitions of four objects into two unordered pairs. The sheet and seam-color labels therefore admit a $\Pi_3$-equivariant identification with the affine points and directions above.
\end{theorem}
\begin{proof}
Two affine points determine a unique nonzero difference and the corresponding line. An affine point and a direction determine one affine line, and two direction points determine the line at infinity. This gives the projective plane of order two.

Affine transformations $a\mapsto Ma+b$ preserve the incidence and act on directions by $d\mapsto Md$. There are $4\cdot6=24$ such transformations, acting faithfully on the four affine points, hence as all of $S_4$. Conversely, any incidence automorphism preserving the line at infinity acts faithfully on the affine points and belongs to this same group.

The rotational cube group has order $24$ and acts faithfully on its four body diagonals. Indeed, a linear orthogonal map preserving all four diagonal lines acts by signs on representatives; the linear relation among any four such representatives has all coefficients nonzero, forcing the signs to agree. The map is therefore $I$ or $-I$, and only $I$ preserves orientation. This gives $S_4$, with kernel $\{\pm I\}$ for the full cube group. For the remaining labels, represent a diagonal by a sign triple modulo overall sign. For axis $i$, the sign of the product of the other two coordinates partitions the four diagonals into two pairs. Signed permutations transport these three partitions exactly as they transport axis classes. These are the three parallel classes of the affine plane.
\end{proof}

This supplies an incidence correspondence. To use the seven labels as imaginary octonion units, fix an oriented Fano multiplication table for $1,e_1,\ldots,e_7$ with $e_i^2=-1$; see~\cite{Baez}.

\begin{proposition}[Criterion for a multiplicative color action]\label{prop:octonionlift}
Suppose a group acts on the seven labels by $\pi_g$, and choose signs $\sigma_{g,i}\in\{\pm1\}$. The maps
$U_g(1)=1$, $U_g(e_i)=\sigma_{g,i}e_{\pi_g(i)}$
define an action by octonion algebra automorphisms precisely when they satisfy the group law and, for every oriented product $e_ie_j=\varepsilon_{ij}e_k$,
\begin{equation}
\sigma_{g,i}\sigma_{g,j}e_{\pi_g(i)}e_{\pi_g(j)}
=\varepsilon_{ij}\sigma_{g,k}e_{\pi_g(k)}.
\end{equation}
\end{proposition}
\begin{proof}
Multiplicativity implies the identities on basis products. Conversely, they, the unit relation, and $e_i^2=-1$ preserve every basis product. Bilinearity gives preservation of every product, and signed permutation maps are invertible. The group law gives the claimed action.
\end{proof}

Geometric symmetry determines how colored strata are permuted. Fano incidence determines the triples involved in multiplication, and the oriented product identities determine the admissible signed lifts. These are successive enrichments. The finite group $\Pi_3$ and the nonassociative algebra $\OO$ are different types of objects; their relevant relationship is such an action. For the geometry, arithmetic, and symmetry of quaternion and octonion algebras, see \cite{ConwaySmith2003}.
The structural theory of octonions and their relationship to
exceptional groups is developed in \cite{SpringerVeldkamp2000}.
In the present construction, the multiplication law is additional
algebraic data whose compatibility with the geometric action must
be checked.

\section{Using the framework as a reference}\label{sec:referenceuse}
\purpose{Collect the objectives in a usable construction order and identify the output of each stage.}

A realization starts with a seed and its attaching maps. It specifies a metric, an operator domain, and any local systems or color data. Local states are selected in that analytic setting, and their seam traces determine whether they define a global state. An appropriate theta or Maass transform gives the Mellin completion. A tower consists of these realizations together with maps transporting their chosen structure.

\begin{center}
\renewcommand{\arraystretch}{1.2}
\begin{tabular}{@{}p{0.19\textwidth}p{0.35\textwidth}p{0.35\textwidth}@{}}
\toprule
Objective & Specified structure & Mathematical output\\
\midrule
Incidence & Sheets, attaching maps, symmetry & An enriched seed and its marked channels\\
Propagation & Metric, operator, boundary data & Local spectral states\\
Compatibility & Common traces and flux law & Global state or seam matching defect\\
Completion & Kernel and Mellin convention & Completed state and its divisor\\
Change of rank & Bonding maps and state transitions & Inverse system of spectral states\\
Duality & Compatible Fourier maps and descent & Reflected universal completed state\\
\bottomrule
\end{tabular}
\end{center}

The invariant Kirchhoff seam produces $\theta_\Gamma$ and $\Lam$. The sheet energy supplies harmonic extension and transmission. The cusp equation supplies $K_{ir}$ propagation; modular invariance supplies the Maass period reflection. The reciprocal-pair theorem gives the higher-rank reflection, and the rank-two residue calculation gives a model transition. The assembly theorem states how these mechanisms become compatible over $\LL$.

Every arrow has a domain, a normalization, and a compatibility law. This permits creative geometric realizations while preserving an exact account of their conclusions. The three-seed is the distinguished meeting point: a finite octahedral incidence pattern carrying local waves, an explicitly solvable invariant theta channel, and the rank-one completed zeta state through which the two sides of the tower are joined.

\clearpage
\appendix
\section{Notation and operator conventions}\label{app:notation}
\renewcommand{\arraystretch}{1.13}
\begin{longtable}{@{}p{0.23\textwidth}p{0.70\textwidth}@{}}
\toprule Symbol & Meaning\\\midrule\endfirsthead
\toprule Symbol & Meaning\\\midrule\endhead
$X,V,\cP$ & Boundary cube, its vertices, and four opposite-vertex pairs.\\
$C_{\partial X}$ & Six face crosses, forming a subdivided octahedral graph.\\
$\mathfrak S_n$ & Enriched seed at level $n$.\\
$\cI_n,\Gamma_n,\Pi_n$ & Sheeted realization, seam, and symmetry group.\\
$\Pi_3$ & Full octahedral group $O_h\cong B_3$, of order $48$.\\
$\FS(X)$ & Full foliated incidence system containing the distinguished seed.\\
$\ell$ & Common length in the equilateral seam model.\\
$A_\Gamma$ & Positive Kirchhoff Laplacian on the metric seam.\\
$A_\Gamma^\G$ & Its restriction to invariant scalar functions.\\
$\mathsf H_\Gamma$ & Normalized operator $\ell^2A_\Gamma^\G/(4\pi)$, with eigenvalues $\pi m^2$.\\
$Z_\Gamma,\theta_\Gamma$ & Invariant heat trace and reconstructed kernel $2Z_\Gamma-1$.\\
$A_3,q_3$ & Sheet operator and its closed energy form.\\
$\cD_\Gamma(\lambda)$ & Seam Dirichlet-to-Neumann operator.\\
$\eta=1/2+ir$ & Leafwise principal parameter; eigenvalue $\eta(1-\eta)$.\\
$\beta=1+i\varrho$ & Bulk principal parameter; eigenvalue $\beta(2-\beta)$.\\
$s$ & Mellin variable, independent of $\eta$ and $\beta$.\\
$k_n=n-2$ & Theta rank at level $n$.\\
$w=s/k$ & Common reflection coordinate at theta weight $k$.\\
$\Lam(s)$ & Meromorphic completion $\pi^{-s/2}\Gamma(s/2)\zeta(s)$.\\
$\xi(s)$ & Entire completion $s(s-1)\Lam(s)/2$.\\
$\Theta,\cZ$ & Families of theta data and completed spectral data.\\
$\widehat{\mathfrak S}_\pm$ & Two inverse-limit towers.\\
$\LL,D$ & Doubled geometric base and its involution.\\
$T_n,U_n$ & Theta and completed-state transitions.\\
$J_n,R_n$ & Theta duality and completed-state reflection.\\
\bottomrule
\end{longtable}

All Laplacians are positive: $A=-\Delta$ for $\Delta=\operatorname{div}\nabla$. Inner products are linear in their first argument. Matching seam traces are part of the form domain. A prime on a spectral trace omits the zero eigenspace. Euler's gamma function is $\Gamma(z)$; the geometric seams carry subscripts $\Gamma_n$.

\section{Reference formulas and convergent domains}\label{app:formulas}
\begin{align}
\vartheta(t)&=\sum_{m\in\Z}e^{-\pi m^2t},&&t>0,\\
\vartheta(t)&=t^{-1/2}\vartheta(1/t),&&t>0,\\
\Lam(s)&=\frac12\int_0^\infty(\vartheta(t)-1)t^{s/2}\frac{\dd t}{t},&&\Re s>1,\\
\Lam(s)&=\pi^{-s/2}\Gamma(s/2)\zeta(s),&&\Re s>1.
\end{align}
For $r\in\R$, $n\ge1$, and $\Re s>0$,
\begin{equation}
\int_0^\infty K_{ir}(2\pi ny)y^{s-1}\dd y
=\frac{\pi^{-s}n^{-s}}4\Gamma\!\left(\frac{s+ir}{2}\right)
\Gamma\!\left(\frac{s-ir}{2}\right).
\end{equation}
Formula~\eqref{eq:splitMellin} extends the zeta identity meromorphically. A Fourier series can be integrated termwise where coefficient growth and kernel estimates give absolute convergence, or after a separately justified continuation. For a reciprocal pair of weight $k$, the initial domain is $\Re s>k$ and the continuation is~\eqref{eq:pairedMellin}. Product kernels have product half-plane domains. These classical identities may also be located in~\cite{DLMFTheta,DLMFZeta}.

\section{Logical dependencies}\label{app:dependencies}

The diagram records dependencies between constructions, rather than maps between their underlying geometric spaces.
\begin{center}
\begin{tikzpicture}[node distance=0.9cm and 0.55cm,
 box/.style={draw=blue!35!black,rounded corners=2pt,align=center,inner sep=7pt,text width=4.0cm,font=\small},
 >={Latex[length=2mm]}]
\node[box] (seed) {Seed incidence and\\ equivariant enrichment};
\node[box,below left=of seed] (graph) {Invariant seam operator\\ Theorem~\ref{thm:interval}};
\node[box,below right=of seed] (sheets) {Sheet energy and matching\\ Theorem~\ref{thm:sheetoperator}};
\node[box,below=of graph] (theta) {Theta and completed zeta\\ Theorem~\ref{thm:seamzeta}};
\node[box,below=of sheets] (cusp) {Cusp and automorphic analysis\\ Theorems~\ref{thm:cusp}, \ref{thm:maassFE}};
\node[box,below=4.8cm of seed] (trans) {Specified transitions\\ and duality compatibility};
\node[box,below=of trans] (assembly) {Universal completed state\\ Theorem~\ref{thm:assembly}};
\draw[->] (seed)--(graph);\draw[->] (seed)--(sheets);
\draw[->] (graph)--(theta);\draw[->] (sheets)--(cusp);
\draw[->] (theta.south) |- (trans.west);\draw[->] (cusp.south) |- (trans.east);
\draw[->] (trans)--(assembly);
\end{tikzpicture}
\end{center}

Theorem~\ref{thm:ranktwo} supplies an explicit transition at the penultimate stage. Section~\ref{sec:divisors} applies to completed states with a specified meromorphic realization. The local-system and Fano structures of Section~\ref{sec:enrichment} enrich the initial datum and can change the operator and state spaces.

\end{document}